\documentclass[11pt]{amsart}

\usepackage{palatino} 
\usepackage{hyperref}
\usepackage{tikz}
\usetikzlibrary{positioning}

\usepackage{amsmath,amssymb,amsthm,verbatim, amsfonts,amscd,flafter,epsf, epsfig,graphicx,verbatim,pinlabel,mathrsfs}
\usepackage[all]{xy}
\usepackage{epsf}
\usepackage[abs]{overpic}
\usepackage{epstopdf}

\usepackage{xcolor}
\usepackage{mathtools}
\usepackage{scalerel}

\definecolor{bred}{rgb}{1,0,.2}
\definecolor{blue}{rgb}{0,0,1}

\newtheorem{theorem}{Theorem}[section]
\newtheorem{lemma}[theorem]{Lemma}

\newtheorem{corollary}[theorem]{Corollary}
\newtheorem{proposition}[theorem]{Proposition}

\theoremstyle{definition}
\newtheorem{definition}[theorem]{Definition}

\newtheorem{remark}[theorem]{Remark}

\def\leave#1{{}}

\def\d{\mathit{d}}

\def\ms#1{{\textcolor{olive}{#1}}}

\def\H{\mathcal{H}}

\def\tb{{\mathit{tb}}}

\title{The Giroux Correspondence in dimension 3}
\author[Joan Licata]{Joan Licata}
\address{Mathematical Sciences Institute, Australian National University}
\email{joan.licata@anu.edu.au}

\author[Matthias Scharitzer]{Matthias Scharitzer}
\address{Uppsala University}
\email{matthias.scharitzer@math.uu.se}

\author[Vera V\'ertesi]{Vera V\'ertesi}
\address{University of Vienna}
\email{vera.vertesi@univie.ac.at}

\begin{document}

\maketitle


\begin{abstract} 
This paper proves the Giroux Correspondence in dimension three using Heegaard splittings of contact manifolds.  In \cite{LV}, the authors proved the Giroux Correspondence for tight contact $3$-manifolds via convex Heegaard surfaces, and simultaneously, \cite{BHH} gave an all-dimensions proof of the Giroux Correspondence by generalising convex surface theory to higher dimensions.  This paper extends the Heegaard splitting approach to arbitrary (not necessarily tight) contact 3-manifolds in order to provide a proof accessible to a low-dimensional audience.  The proof assumes a result announced in \cite{BHH} that classifies moves relating  bypass decompositions for  isotopic contact structures on cobordisms that are topological products; in the Appendix, we prove this result in the 3-dimensional setting.

\end{abstract}

\section{Introduction}
Topology and contact geometry enjoy a remarkable relationship in dimension $3$. Contact structures are defined in the smooth category with differential geometric language, and yet they mirror fundamental topological properties of the underlying $3$-manifold. An open book decomposition is arguably the most striking incarnation of this relationship: up to diffeomorphism, a contact manifold  is determined simply by a surface mapping class.  For the last twenty years, 3-dimensional contact geometry has been shaped by the equivalence between contact structures and open books first proposed by Giroux in the early 2000s:

\begin{theorem}\label{thm:GC}[Giroux Correspondence \cite{Giob}] Two open book decompositions support isotopic contact structures if and only if they  are related by a sequence of positive (de)stabilisations.
\end{theorem}

Because topology and contact geometry are so closely entwined, it's natural to probe whether topological techniques can be squeezed to produce meaningful contact geometric outputs.  In this spirit, in 2023 the authors developed new tools to study $3$-dimensional contact manifolds via their Heegaard decompositions. Using this enhanced relationship between contact structures and Heegaard splittings, we proved the Giroux Correspondence for all tight contact $3$-manifolds \cite{LV}. Simultaneously, Breen-Honda-Huang used generalised convex surface theory to prove the Giroux Correspondence for all contact structures and in all dimensions \cite{BHH}.  While the work of Breen-Honda-Huang fully covers the $3$-dimensional case, the generality of their techniques lends their proof a different flavour.

This paper provides a complete proof of the Giroux Correspondence in three dimensions, further developing the Heegaard splitting approach to this topic and using an essential result of \cite{BHH} (Theorem~\ref{thm:bypass1p} below) to bridge the gap to overtwisted manifolds. Not only will this proof be accessible to a low-dimensional audience, but it also illustrates the effectiveness of applying convex surface techniques in concert with Heegaard splittings, an approach which we hope will see more applications in the future. At the request of the referee, we have also included a proof of Theorem~\ref{thm:bypass1p} in Appendix~\ref{Appendix}.  Although similar in spirit to the argument in \cite{BHH}, some techniques specific to low dimensions affect the flavour of the argument.

\subsection{Heegaard splittings and open book decompositions}
Any open book decomposition of a manifold gives rise to a canonical Heegaard decomposition.  Heegaard splittings produced thus are called \textit{convex Heegaard splittings}, and one may directly recover an open book supporting $\xi$ from any convex Heegaard splitting of $(M, \xi)$.  The innovation in \cite{LV} is a procedure called \textit{refinement} that produces an open book from a more general class of Heegaard splittings of a contact manifold.  Although this open book decomposition is not uniquely defined, any two open books produced thus are related by a sequence of positive (de)stabilisations -- exactly the equivalence appearing in the Giroux Correspondence.  It follows that the Giroux Correspondence can be rephrased as a statement about convex Heegaard splittings (Theorem~\ref{thm:mainG}). Our proof of the Giroux Correspondence for tight contact manifolds relied on showing that many operations  one might naturally perform on a Heegaard splitting of a contact manifold in fact preserve the positive stabilisation class of the associated open book.

As explained next, these techniques were insufficient in the case of overtwisted manifolds.  In order to assign an open book decomposition to a contact manifold with a Heegaard splitting, we required both handlebodies in the Heegaard splitting to be tight.  When two convex Heegaard surfaces are isotopic, there is a sequence of convex Heegaard surfaces interpolating between them.  Our approach tried to assign a sequence of open book decompositions to this sequence of Heegaard splittings.  However, even when both original Heegaard splittings have tight handlebodies,  intermediate Heegaard splittings in this sequence may nevertheless have overtwisted handlebodies.   Since we could not preclude this for arbitrary contact $3$-manifolds, our proof of the Giroux Correspondence was limited to tight manifolds.

The present work bypasses this limitation by considering highly stabilised Heegaard splittings where the handlebodies are always tight.  We introduce the term ``bridging'' for a particular method of constructing such a splitting, and we show that different choices made in constructing a bridge splitting preserve the positive stabilisation class of the associated open book.  Using bridge splittings,  we prove that any pair of convex Heegaard splittings for a fixed contact manifold correspond to open books that are related by a sequence of positive (de)stabilisations.

\subsection{A note on notation} 
The arguments in this paper involve sequences of Heegaard splittings of varying types, each of which is denoted by some sort of  ``$\H$''.  We have attempted to define notation locally and maintain these choices consistently, but we briefly outline the system here for reference.  Every Heegaard surface $\Sigma$ is required to be convex; Heegaard splittings subject to no further restrictions are denoted by $\H=(\Sigma, U, V)$, where $U$ and $V$ are the handlebodies, oriented so that $\partial U=\Sigma$.  Elements of a pair of Heegaard splittings are distinguished by the addition of a prime: $\H=(\Sigma, U, V)$ and $\H'=(\Sigma', U', V')$.  Following the notation in \cite{LV}, the refinement of $\H$ is denoted by $\widetilde{\H}=(\widetilde{\Sigma},\widetilde{U},\widetilde{V})$.

Starting in Section~\ref{sec:hsbyp}, we will stabilise Heegaard splittings by attaching certain contact $1$-handles found in the complement of fixed handlebodies.  This process is called \textit{bridging}, and we denote the bridge of $\H$ by $\widehat{\H}=(\widehat{\Sigma},\widehat{U},\widehat{V})$.  When we  bridge only a single handlebody of the original splitting, we denote the result by $\widehat{\H}^U$ or $\widehat{\H}^V$, as appropriate.  Distinct bridged splittings may be distinguished by primes: $\widehat{\H}$ and $\widehat{\H}'$. 

\subsection*{Acknowledgment} This research was supported in part by the Austrian Science Fund (FWF) P 34318. MS was supported by the Swedish Research Council, VR 2022-06593, Centre of Excellence in Geometry and Physics at Uppsala University and VR 2024-04417, project grant. For open access purposes, the author has applied a CC BY public copyright license to any author-accepted manuscript version arising from this submission. The authors would also like to thank  Joseph Breen, Nataliya Goncharuk and Ko Honda for helpful discussions on the topic.

\section{Background}\label{sec:prelim} We assume the reader is familiar with basic notions in contact geometry, but we use this section to collect some tools and terminology that will be essential for the rest of the paper. Section~\ref{sec:convex} provides a quick review of convex surface theory and introduces \textit{weak contactomorphism} (Definition~\ref{def:wc}) as an equivalence relation on contact manifolds with convex boundary.  Section~\ref{sec:ob} describes the relationship between open book decompositions and convex Heegaard splittings (Definition~\ref{def:convexhs}) of a contact manifold. Section~\ref{sec:bypass} presents the contact handle model for bypass attachment.  Finally, Section~\ref{sec:bypassdecomp} states the key theorem  which classifies distinct bypass decompositions of a product cobordism. This statement will be proved in the  \ref{Appendix}. Also, Lemma~\ref{lem:rotex} gives a new perspective on the familiar operation of bypass rotation.

\subsection{Legendrian Graphs} In Section~\ref{sec:splittingsandgc} will consider contact Heegaard decompositions constructed as neighbourhoods of Legendrian graphs. Working with these Legendrian graphs requires a 1-parameter version of the Legendrian approximation theorem that was originally proved for Legendrian knots in the standard contact structure by Fuchs and Tabachnikov \cite{FuchTab}.

\begin{theorem}\label{thm:Legapprox1} For $t\in[0,1]$, let $G_t$ be a 1-parameter family of spacial graphs in a contact manifold $(M,\xi)$.  Then $G_t$ can be $C^0$-approximated by a 1-parameter family of Legendrian graphs. 
\end{theorem}

That is, for any smooth  choice of neighbourhoods $N_t$ of $G_t$,  we can choose a 1-parameter family of Legendrian graphs $L_t$ such that $L_t$ is smoothly isotopic to $G_t$ inside $N_t$.

Note that at a vertex $v$ of a Legendrian graph,  all the edges must be tangent to the  contact plane $\xi_v$.  It follows that the edges have a cyclic order that must be preserved throughout a Legendrian isotopy.  This explains the final hypothesis in the next result.  

\begin{corollary}\label{thm:Legapprox2} Suppose the Legendrian graphs $L_0$ and $L_1$ are smoothly isotopic, and suppose, too, that the cyclic orders of the edges agree at every vertex. Then after finitely many positive and negative stabilistions, $L_0$ and $L_1$ become Legendrian isotopic.
\end{corollary}

\begin{proof}[Proof of Theorem \ref{thm:Legapprox1}]
At each vertex of $G_t$, fix a cylic order on the incoming edges.  Isotope $G_t$ inside $N_t$ so that near each vertex $v$, the edges are Legendrian in a Darboux-ball $D_v(t)$ and realise the given cyclic order for all $t\in[0,1]$. Now extend the isotopy of $G_t$ to an ambient isotopy $\Phi_t$ of $M$ that brings the Daurboux balls $D_v(0)$ to $D_v(t)$ via a contactomorphism. Pull the contact structure back by the ambient isotopy. This produces a fixed graph $G=\Phi^{-1}_t(G_t)$ in a 1-parameter family of contact structures $\xi_t=(\Phi_t)_*\xi$, all of which agree near the vertices of $G$. Choose a neighbourhood  $N$ of $G$ contained in the intersection of the neighbourhoods $\Phi_t^{-1}(N_t)$.

Now restrict to an edge $e=(u,v)$ and follow the usual procedure for Legendrian approximation, as follows. Using standard transversality arguments, a $C^\infty$-small isotopy of $\xi_t$ ensures that away from $\coprod D_v$,  the edge $e$ is transverse to $\xi_t$ for all but  finitely many times $0<t_1<\dots<t_n<1$.  Furthermore, at each such $t_j$,   exactly at one point in the interior of $e$ has a non-degenerate tangency with $\xi_{t_j}$.

In  time intervals away from the singular $t_j$ values, take standard neighbourhoods contactomorphic to $(I \times D^2,\ker(ds+r^2d\varphi))$, allowing the radius to vary in $t$. By the compactness of each interval $[t_j+\varepsilon, t_{j+1}-\varepsilon]$, there exists a radius $R$ defining a neighbouhood contained in $N$ for all $t$.  Use a leaf of the characteristic foliation on $T_R=I\times \{r=R\}$ together with a Legendrian model in $D_u\cup D_v$ to approximate $e$ in $N$. Different choices of $R$ correspond to positive (de)stabilisations of the Legendrian approximation. 

At a singular times $t_j$, the interior of the edge $e$ is tangent to $\xi_{t_j}$ at a point $p$.  That is,  the edges $\Phi_{t-\varepsilon}(e)$ and $\Phi_{t+\varepsilon}(e)$ are related to each other by transverse stabilisation or destabilisation, depending on the sign of $\frac \d{d t} \alpha_{t}(T_pe)\vert_{t=t_j}$. 
For simplicity, assume $\Phi_{t+\varepsilon}(e)$ is a transverse stabilistion of $\Phi_{t-\varepsilon}(e)$. Then a standard neighbourhood of $\Phi_{t+\varepsilon}(e)$ is contained in the standard neighbourhood of $\Phi_{t-\varepsilon}(e)$, and their Legendrian approximations as above are related to each other via a single negative Legendrian stabilisation and some positive Legendrian (de)stabilistions that depend on the $R$-values chosen for $[t_{j-1}+\varepsilon,t_j-\varepsilon]$ and $[t_j+\varepsilon,t_{j+1}-\varepsilon]$. 


Thus we have proved that Legendrian isotopy and positive and negative Legendrian stabilisations suffice to approximate $G_t$. {Corollary~\ref{thm:Legapprox2}}  follows from performing the stabilisations before the isotopies. 
\end{proof}

\subsection{Convex surfaces}\label{sec:convex}
The reader looking for a more thorough introduction to convex surface theory is referred to \cite{Ma} or the lecture notes \cite{Econ,Hcon}. 

Recall that a vector field $\mathcal{V}$ in a contact manifold $(M, \xi)$ is \textit{contact} if the flow of the vector field preserves $\xi$.  A surface is \textit{convex} when there exists a  contact vector field $\mathcal{V}$ transverse to $\Sigma$. 
The essential feature of a convex surface $\Sigma$ is that it provides a combinatorial characterisation of the contact structure in an $I$-invariant neighbourhood $\nu(\Sigma)$; the product structure on $\nu(\Sigma)$ is provided by the flow of $\mathcal{V}$.  The \textit{dividing set} on $\Sigma$ is a separating multicurve $\Gamma_\Sigma$ consisting of the points where $\mathcal{V}$ lies in $\xi$.  Changing the vector field that certifies $\Sigma$ as convex changes $\Gamma_\Sigma$ only by isotopy, and any isotopy keeping  $\Sigma$ convex also preserves $\Gamma_\Sigma$ up to isotopy.  Furthermore, any $\Gamma'$ isotopic to $\Gamma_\Sigma$ may be realised as a dividing set after performing an isotopy of $\Sigma$ within an $I$-invariant neighbourhood.  

Convex surface theory is particularly useful for studying contact manifolds with boundary, and we will assume throughout that  the boundary of $(M, \xi)$ is convex.  Since the dividing set characterises  $\xi$ near $\partial M$,  we would like to be able to glue $(M, \xi, \Gamma_{\partial M})$ to $(M',  \xi', \Gamma_{\partial M'})$ whenever there is a diffeomorphism taking $(\partial M, \Gamma_{\partial M})$ to $(\partial M', \Gamma_{\partial M'})$.   Because we are working in the smooth category, gluing is achieved by identifying  contact neighbourhoods of the boundaries; for later convenience, we introduce the following equivalence relation on contact manifolds with boundary:

\begin{definition}\cite[Definition 2.3.]{LV}\label{def:wc} 
Two contact manifolds $(M,\xi)$ and $(M',\xi')$ with  convex boundary are \emph{weakly contactomorphic} if there is a contact embedding $\iota\colon M\hookrightarrow M'$  such that $\iota(\partial M)$ is convexly isotopic to $\partial M'$. \end{definition}

\begin{remark} We note that if $\iota:M\rightarrow M'$ is a contact embeddeding as above, then $M\setminus \iota({M'})$ is a weak identity morphism in the sense of \cite{BT}. 
\end{remark}

Henceforth, we will be interested in contact manifolds with boundary only up to weak contactomorphism.  Two (weak contactomorphism classes of) manifolds may be glued along diffeomorphic boundary components if and only if  there is a  diffeomorphism respecting the dividing sets. Gluing is well defined as an operation on weak contactomorphism classes of manifolds; see \cite[Proposition 2.8.]{LV}.

As a special case of weakly contactomorphic manifolds, we also introduce the following:

\begin{definition}\cite[Definition 2.4.]{LV} Two embedded codimension-$0$ submanifolds $N_1, N_2\subset (M, \xi)$ with convex boundary are \textit{weakly (contact) isotopic} if there is an isotopy $N_s$ between them such that $\partial N_s$ remains convex throughout.
\end{definition}

\subsection{Open book decompositions of contact manifolds}\label{sec:ob}

Let $(B,\pi)$ denote an open book decomposition of $M$.  That is, $B$ is an oriented link and $\pi:M\setminus B\rightarrow S^1$ is a fibration of the complement of $B$ by Seifert surfaces for $B$. When $M$ is equipped with a contact structure $\xi$, we assume any open book decomposition supports $\xi$ in the sense of \cite{TW}.

An open book decomposition $(B,\pi)$ induces a Heegaard splitting $(\Sigma, U,V)$ where each handlebody consists of half the pages: $U=\overline{\pi^{-1}[0,1/2]}$ and $V=\overline{\pi^{-1}[1/2,1]}$. The Heegaard surface \[\Sigma=\overline{\pi^{-1}(0)}\cup-\overline{\pi^{-1}(1/2)}\] 
is naturally convex with dividing curve $\Gamma_\Sigma=B=\partial \big(\pi^{-1}(0)\big)$.  An arc properly embedded on a page of the open book produces  a properly embedded disc in the corresponding handlebody, and in fact, such discs characterise the Heegaard splittings induced by open books. 

More precisely, a disc $A$ properly embedded in a contact manifold with convex boundary $(\Sigma, \Gamma_\Sigma)$ is a \textit{product disc} if $|\partial A \cap \Gamma_\Sigma|=2$.  A \textit{disc system} for a handlebody is a set of properly embedded discs that cut it into a ball. Torisu introduced the notion of a convex Heegaard splitting \cite{Torisu}:

\begin{definition}\label{def:convexhs} A Heegaard splitting $(\Sigma,U,V)$ of $(M,\xi)$ is   \emph{convex} if 
\begin{itemize}
    \item $\Sigma$ is convex with dividing curve $\Gamma_\Sigma$;
        \item $\xi\vert_U$ and $\xi\vert_V$ are both tight;  and
    \item each of the manifolds $U$ and $V$ admits a  system of product discs.
    \end{itemize}
\end{definition}

The discs in a convex Heegaard splitting determine a product structure  on each handlebody, and hence,  an open book decomposition with binding $\Gamma_\Sigma$. Moreover, this relation is one-to-one up to isotopy. 

\begin{proposition}
    \label{lem:obvshs}\cite[Proposition 3.4.]{LV} 
    Let $(M, \xi)$ be a contact manifold and let $(B,\pi)$ and $(B',\pi')$
be open books supporting $\xi$. Then the induced convex Heegaard splittings are isotopic via an
isotopy  keeping the Heegaard surface convex if and only if $(B,\pi)$
and $(B',\pi')$ are isotopic through a path of open books supporting $\xi$.
\end{proposition}

Because  a convex Heegaard splitting  determines an open book decomposition, a convex splitting allows one to reconstruct the contact structure.  Specifically, suppose $(\Sigma, \Gamma_\Sigma, \{\alpha_i\}, \{\beta_i\})$ is a convex Heegaard diagram for $(M, \xi)$ where the meridional circles bound product discs in the respective handlebodies and the restriction of $\xi$ to each handlebody is known to be tight.  Then  up to contact isotopy, there is a unique contact structure structure compatible with this data.

\subsection{Bypass attachment}\label{sec:bypass}
Next, we turn to an operation that changes the weak contactomorphism class of a contact manifold while preserving the topological type. Adopting the perspective of \cite{HKM, Oz}, we  define bypass attachment  as an operation performed by gluing contact handles.  

As in the topological setting, a contact $k$-handle is a $D^k\times D^{3-k}$ glued along $\partial D^k \times D^{3-k}$ to the boundary of an existing manifold $M$.  Additionally, each $k$-handle is required to smooth to a tight $3$-ball, and the dividing set on the attaching region is prescribed. A $1$-handle attachment is specified by an $S^0$ on the dividing set $\Gamma_{\partial M}$, while a $2$-handle attachment is specified by an $S^1$ that intersects $\Gamma_{\partial M}$ twice transversely. Local models for handle attachment are well established in the literature; see, for example, \cite{Gi} or \cite{Oz}.  Although each handle is a cornered manifold, the gluing process is defined so that the manifold is smooth after each attachment.

Let $(\Sigma, \Gamma_\Sigma)$ be a convex surface. A Legendrian arc $c$ embedded in  $\Sigma$ is a \textit{bypass arc} if $\partial c$ lies on $\Gamma_\Sigma$ and the interior of $c$ intersects $\Gamma_\Sigma$ once transversely. As explained next, we may associate a new contact manifold with convex boundary to any bypass arc on the convex boundary of a contact manifold.  

 \begin{figure}[h]
\begin{center}
\includegraphics[scale=1.2]{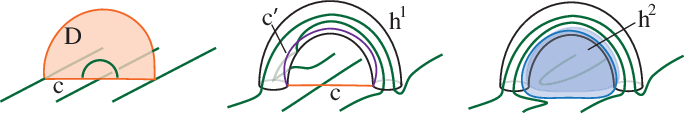}
\caption{ Left: a bypass half-disc attached along the arc $c$.  Centre: the attached $1$-handle is a regular neighbourhood of an arc Legendrian isotopic to $c$.  Right: attaching the shaded $2$-handle along $c\cup c'$ completes the bypass attachment. }\label{fig:bypassmodel}
\end{center}
\end{figure}
 
\begin{figure}[h]
\begin{center}
\includegraphics[scale=1.5]{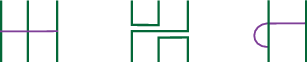}
\caption{Left: Local model for a bypass attachment arc.  Centre: The new dividing set after attaching a bypass along the arc shown on the left.  Right: A trivial bypass  arc.  Attaching a bypass along the arc shown preserves the dividing set up to isotopy. }\label{fig:bypass}
\end{center}
\end{figure}

The steps described here are illustrated in Figure~\ref{fig:bypassmodel}. First, attach a minimally twisting contact $1$-handle $h^1$ along the $0$-sphere $\partial c$.  Slightly abusing notation, continue to denote the remaining portion of the bypass arc by $c$.  Let $c'\subset \partial h^1$ be a Legendrian arc with endpoints on $c$ and such that $\tb(c\cup c')=-1$.  Attach a contact $2$-handle $h^2$ along $c \cup c'$.    Topologically, the handles cancel, but in general the dividing set on the new boundary is not isotopic to $\Gamma_\Sigma$.  See Figure~\ref{fig:bypass}.  We often write $X$ to denote the pair of handles $h^1\cup h^2$; $X$ is itself a neighbourhood of a \textit{bypass half-disc} $D$ whose cornered Legendrian boundary consists of the bypass arc $c$ and a Legendrian core of $h^1$. Coordinate models for this construction are provided in \cite{Oz}, and the dividing sets at each stage of the construction are shown in Figure~\ref{fig:bypassmodel}.    

Suppose that $(M, \xi, \Gamma_{\partial M})$ is a contact manifold with convex boundary and let $c$ be a bypass arc on $\partial M$.  Then the weak contactomorphism class of $M \cup (h^1\cup h^2)$ is preserved under isotopies of $c$ through bypass arcs.

A bypass is \textit{trivial} when the bypass arc $c$ on $\Sigma$ matches the one shown on the right in Figure~\ref{fig:bypass}.  That is, the arc $c$ turns left after the interior crossing with $\Gamma_\Sigma$ to cobound a bigon with the same arc of $\Gamma_\Sigma$.

 Attaching a trivial bypass preserves the convex isotopy class of the surface.  If $\Sigma$ bounds a contact submanifold $(N, \xi)$, then attaching a trivial bypass also preserves the weak contact isotopy class of $N$. {As shown in \cite[Lemma~5.1.]{HKMrv}, when the bypass arc $c$ is trivial on $\Sigma\subset (M, \xi)$, the corresponding bypass half-disc $D$ always exists as a submanifold of $M$.

 \subsection{Bypass decompositions}\label{sec:bypassdecomp} Bypasses allow us to classify distinct contact structures on a fixed topological manifold.  Adopting cobordism notation, we write $\partial_+M$ (respectively,  $\partial_-M$) to distinguish the outward-oriented (respectively, inward-oriented)  boundary of a topological $M$. 
 
 \begin{theorem}\label{thm:bypass}\cite[Section 3.2.3]{Hgluing} Let $(\Sigma \times I, \xi)$ be a topological product with convex boundary.  Then there exists a sequence of bypasses $X_i:=h_i^1\cup h_i^2$ such that $(\Sigma \times I, \xi)$ is weakly contactomorphic to
\[ \nu(\Sigma)\cup \bigcup_i X_i,\]
where $\nu(\Sigma)$ is an $I$-invariant half-neighbourhood of $\Sigma \times \{-1\}$ and the bypass arc $c_j$ 
for $X_j=h_j^1\cup h_j^2$ lies on 
\[ \partial_+\big( \nu(\Sigma)\cup \bigcup_{i=1}^{j-1} X_i\big).\]
\end{theorem}

Theorem \ref{thm:bypass} has an alternative presentation as a decomposition theorem.  Given a contact manifold $(\Sigma \times I, \xi)$ with convex boundary, consider an embedding 
\[ \iota:\nu(\Sigma)\cup \bigcup_i X_i\hookrightarrow (\Sigma \times I, \xi)\]
as in Definition~\ref{def:wc}.
It follows that 
contact handles $h_i^1$ and $ h_i^2$ may be identified as submanifolds of the original $(\Sigma \times I, \xi)$.  This identifies a sequence of bypass arcs $c_j$ on the embedded convex surfaces  $\Sigma_{j-1}:=\iota(\partial_+\big( \nu(\Sigma)\cup \bigcup_{i=1}^{j-1} X_i\big)\big)$. 
 Similarly, there exist embedded bypass half-discs $D_i$ in $\Sigma \times I$.  The surface $\Sigma_j$ is constructed by isotoping  $\Sigma_{j-1}$ across a neighbourhood of $D_j$.  
 
 We write $(\Sigma\times I, \xi)=X_1 \ast X_2\ast\dots \ast X_f$ to denote this decomposition of $(\Sigma \times I, \xi)$. 

The factorisation provided by this theorem  is not unique.  As mentioned before, one can isotope any bypass attaching arc through bypass arcs.  In addition, there are  two further moves that change the decomposition while preserving the contact manifold: Trivial Insertion and Far Commutation. 

\noindent\textit{[TI]: Trivial Insertion} 
Given a bypass decomposition $X_1 \ast \dots \ast X_f$, let $c_T$ denote a trivial bypass arc on $\Sigma_{j}$.  Then a bypass attachment along $c_T$ may be inserted into the decomposition:
\begin{center}
\begin{tikzpicture}[node distance = 2cm, thick]%
        \node (1) at (0,0) [above]{$X_1 \ast \dots \ast X_j \ast X_{j+1} \ast\dots \ast X_f$};
       \node (2) at (0,-1)  {$X_1 \ast \dots \ast X_j \ast X_T\ast X_{j+1} \ast\dots \ast X_f$};
 \draw[<->] (1) -- (2);
    \end{tikzpicture}%
    \end{center}

\noindent\textit{[FC]: Far Commutation } Suppose that the bypass arcs $c_j$ and $c_{j+1}$ for $X_j$ and $X_{j+1}$
are disjoint. Then in any  decomposition in which these correspond to consecutive bypass attachments, the order in which these are performed may be exchanged:
\begin{center}
\begin{tikzpicture}[node distance = 2cm, thick]%
        \node (1) at (0,0) [above]{$X_1 \ast \dots \ast X_j \ast X_{j+1} \ast\dots \ast X_f$};
       \node (2) at (0,-1)  {$X_1 \ast \dots \ast X_{j+1} \ast X_{j} \ast\dots \ast X_f$};
 \draw[<->] (1) -- (2);
    \end{tikzpicture}%
    \end{center}

The main theorem of both \cite{BT} and \cite{BHH}  establishes that Trival Insertion and Far Commutation suffice to relate any two  bypass decompositions of topologically trivial products:

\begin{theorem}\cite[Theorem 3.1.2 ]{BHH}\label{thm:bypass1p}
  Any two bypass decompositions of $(\Sigma\times I, \xi)$ are related by Bypass Arc Isotopy and a finite iteration of Trivial Insertion and Far Commutation moves. 
\end{theorem} 

The dimension 3 case was the Main Theorem in \cite{BT}, although presented in different language. For completeness, we give a proof of this statement using independent techniques in the \hyperref[Appendix]{Appendix}.

The term \emph{bypass rotation} is used in several different, but related, senses in the literature.  Here, we consider the phenomenon as presented in \cite{HKMrv} and recall the proof that it is a consequence of the above two moves. Let $c$ be a bypass arc on $\Sigma$. Let $c_T$ be another bypass arc as in Figure~\ref{fig:rot}; i.e.,  subsegments of $c$ and $c_T$ together with two segments of $\Gamma_\Sigma$ cobound a rectangle $R$ that is to the right of $c$ when $c$ is oriented as the boundary of $R$. Let $R'\subset \Sigma$ be a neighbourhood of $R\cup c\cup c_T$ shown shaded in the right hand picture of Figure~\ref{fig:rot}. The arc $c_T$ is said to be a \emph{rotation} of the bypass arc $c$.

 \begin{figure}[h]
\begin{center}
\includegraphics[scale=1.5]{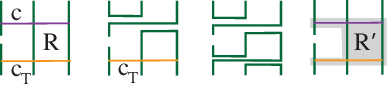}
\caption{ Left: The bypass arc $c_T$ is a rotation of the bypass arc $c$.  Left centre: After attaching a bypass along $c$, $c_T$ becomes a trivial bypass arc. Right centre: The dividing set is preserved up to isotopy by attaching a bypass to $c_T$.  Right: The bypass disc attached along $c_T$ can be constructed in a neighbourhood of 
$R'$ and the bypass disc for $c$.}\label{fig:rot}
\end{center}
\end{figure}

\begin{lemma}\label{lem:rotex}
Suppose that a bypass disc $D$ is attached to the arc $c$ on some surface $\Sigma\subset (M, \xi)$, as shown in Figure~\ref{fig:rot}. Then there also exists a bypass disc $D_T$ attached to $\Sigma$ along $c_T$.  Furthermore, $D_T$ is contained in a neighbourhood of $D$ and $R'$.
\end{lemma}

The second claim of the lemma will be useful in the proof of Proposition~\ref{prop:extindep} below.

\begin{proof} The existence claim follows from a formal convex surface theory argument.  After  a bypass is attached along $c$,  $c_T$ is a trivial bypass arc on the new surface. Since trivial bypasses always exist, it follows that there exists a bypass disc $D_T$ attached along $c_T$.  
The [\textit{TI}] move then implies that $\nu(\Sigma)\cup X_c$ can  be factored as $\nu(\Sigma)\cup (X_c\ast X_{c_T})$. However, the bypass arcs $c$ and $c_T$ are disjoint on $\Sigma$, so the associated bypasses commute.   Applying the [\textit{FC}] move gives the factorisation $\nu(\Sigma)\cup (X_{c_T}\ast X_{c})$. It follows  that the bypass along $c_T$ can be directly attached onto $\Sigma$. 


 The next argument is illustrated in Figure~\ref{fig:existtrivbyp}.  
  \begin{figure}[h]
\begin{center}
\includegraphics[scale=1.2]{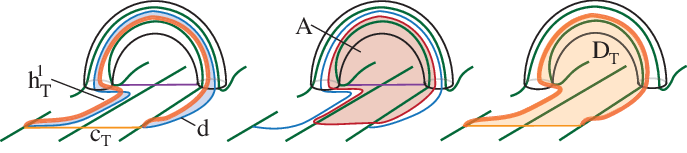}
\caption{ Attach a $1$-handle $h_T^1$ along a push-off of the Legendrian curve $d$.  Then isotope the disc cobounded  by $h_T^1$ and $d$ across the product disc $A$ to produce a bypass disc $D$ cobounded by $h_T^1$ and $c_T$. } 
\label{fig:existtrivbyp}
\end{center}
\end{figure}

To construct $D_T$ directly, consider a Legendrian realisation of the arc labeled $d$ in Figure~\ref{fig:existtrivbyp}.  Perform a Legendrian isotopy so that $d$ is properly embedded in $\nu(\Sigma)\cup X_c$.  Since  $d\cap \Gamma_\Sigma=\partial d$, a regular neighbourhood of this push-off is a contact $1$-handle $h_T^1$ shown in bold on the first figure.  Observe that $h_T^1$ is attached to $\Sigma$ along $\partial c_T$, but the blue disc it cobounds with $d$ is not a bypass disc.  However, we may isotope this disc across the product disc $A$ shown in the centre figure to get a bypass disc $D_T$ cobounded by $h_T^1$ and $c_T$. 
\end{proof}

\section{Heegaard Splittings and the Giroux Correspondence}\label{sec:splittingsandgc}
With useful tools and vocabulary established, we now turn to proving the Giroux Correspondence.  As a first step, Section~\ref{sec:ths} summarises results from \cite{LV} that allow us to restate the theorem about a pair of open books as a theorem about a pair of Heegaard splittings.  Section~\ref{sec:hsbyp} introduces   \textit{bridging} as an operation to stabilise a Heegaard splitting of a contact manifold.  After establishing some properties of bridge splittings, we relate them to refinements in Section~\ref{sec:gc} to prove the main result.

\subsection{Tight Heegaard splittings and refinement}\label{sec:ths} As described in Section~\ref{sec:ob}, each convex Heegaard splitting of a contact manifold corresponds to an open book decomposition  supporting the contact structure.   Here we recall a key idea from  \cite{LV} that extends the class of Heegaard splittings of a contact manifold which determine open books.  

\begin{definition} The Heegaard splitting $(\Sigma, U, V)$ of $(M, \xi)$ is \textit{tight} if $\Sigma$ is convex and the restriction of $\xi$ to each of $U$ and $V$ is tight.
\end{definition}

The definition of a convex Heegaard splitting requires a system of meridional discs each intersecting $\Gamma_\Sigma$ twice.  In contrast, meridional discs in a tight Heegaard splitting may intersect $\Gamma_\Sigma$ an arbitrary number of times.  However, a process called \textit{refinement} allows us to construct a convex splitting from a tight one.

Refinement was introduced in \cite[Section 5.]{LV} . Roughly speaking, refinement stabilises the tight Heegaard splitting by drilling along Legendrian curves embedded in meridional discs.  The drilling process cuts each non-product disc into subdiscs, and the curves are chosen so that these subdiscs are product discs for the stabilised Heegaard splitting $\widetilde{\H}$.   See Figure~\ref{fig:refinement}.  This refinement  $\widetilde{\H}$ is necessarily convex, so via Lemma~\ref{lem:obvshs} one may view refinement as a process which associates an open book decomposition to a tight Heegaard splitting of $(M, \xi)$. 
 \begin{figure}[h]
\begin{center}
\includegraphics[scale=1.2]{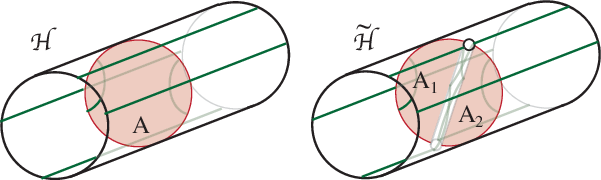}
\caption{ To refine a tight Heegaard splitting $\H$, stabilise it via tunnels drilled through non-product discs to produce a convex splitting $\widetilde{H}$. Here, the meriodonal disc $A$ in $\H$ with $|\partial A \cap \Gamma_\Sigma|=4$  splits into a pair of product discs $A_1$, $A_2$ in $\widetilde{\H}$.} 
\label{fig:refinement}
\end{center}
\end{figure}

There are choices involved in constructing a refinement (e.g., discs, tunnels), but these preserve the positive stabilisation class of the associated open book.

\begin{theorem}\label{thm:refwd}\cite[Lemma 3.8 and Theorem 5.10]{LV} Suppose that $\widetilde{\H}$ and $\widetilde{\H}'$ are refinements of the tight Heegaard splitting $\H$. If $(B, \pi)$ and $(B', \pi')$ are the open book decompositions associated to $\widetilde{\H}$ and $\widetilde{\H}'$, respectively, then $(B, \pi)$ and $(B', \pi')$  are related by a sequence of positive (de)stabilisations.
.
\end{theorem}

In light of Theorem~\ref{thm:refwd}, we  define a \emph{positive stabilisation} of a tight Heegaard splitting as a topological stabilisation of the Heegaard splitting which induces a positive open book stabilisation on the open books associated to the refinements.  

This new language permits us to restate the hard direction of the Giroux Correspondence:

\begin{theorem}[3-dimensional Giroux Correspondence]\label{thm:mainG}
Any two convex Heegaard splittings of $(M, \xi)$ are related by a sequence of positive (de)stabilisations. 
\end{theorem}

Much of the work of \cite{LV} was devoted to showing that various natural operations one may perform on a tight Heegaard splitting are in fact (sequences of) positive stabilisations.   

Positive stabilisations may also be characterised directly.  Choose an arc properly embedded in one of the handlebodies  that is Legendrian isotopic to an arc properly embedded in the closure of $\Sigma\setminus \Gamma_\Sigma$.  Adding a standard contact neighborhood of this arc to the other handlebody is a positive stabilisation of the original Heegaard splitting.  In fact, a related operation has already appeared in the  proof of Lemma~\ref{lem:rotex}.  Suppose that in Lemma~\ref{lem:rotex}  the surface $\Sigma$ is taken to be a Heegaard surface for a tight Heegaard splitting.  In this case,  the $1$-handle $h^1_T$ is a neighbourhood of a Legendrian push-off of $d$, so adding $h_T^1$ to the handlebody bounded by $\Sigma$ is a positive stabilisation of the Heegaard splitting.

The proof of Lemma 3.10 in \cite{LV} shows that in a convex Heegaard splitting, attaching any bypass $1$-handle is a positive stabilisation.  
The following more general statement will be used later.

\begin{theorem}\cite[Lemma 7.1]{LV}\label{thm:lv7.1}
Let $\H=(\Sigma,U,V)$ and $\H'=(\Sigma',U',V')$
be tight Heegaard splittings of $(M,\xi)$ and assume that $U'$ is obtained from $U$ by attaching a single bypass $h^1\cup h^2$.  Let \[\H''=\big(\partial(U\cup h^1), U\cup h^1,M\setminus \text{int}(U\cup h^1)\big)\] be the Heegaard splitting formed by attaching only the $1$-handle associated to this bypass.  
Then the refinements of $\H$, $\H'$ and $\H''$ are all related by sequences of positive (de)stabilisations.

 \end{theorem}

\begin{remark}\label{rmk:extrastab}  Refining a convex Heegaard splitting with respect to a set of product discs preserves the splitting.  However, a convex splitting $\H$ will in general also have non-product meridional discs.  Refining $\H$ with respect to some set of non-product discs will produce a new convex splitting $\widehat{\H}$; then $\H$ and $\widehat{\H}$ admit a common positive stabilisation. 
\end{remark}

\subsection{Bridge Heegaard splittings}\label{sec:hsbyp}

Given a tight Heegaard splitting $\H$, the refinement process described above stabilises $\H$ to produce a convex Heegaard splitting $\widetilde{\H}$.  In this section we introduce an alternative method to turn any smooth Heegaard splitting $\H$ into a convex Heegaard splitting $\widehat{\H}$.  The new Heegaard splitting $\widehat{\H}$ is called the \textit{bridge} of $\H$.

A bouquet of circles is a \textit{skeleton} of a handlebody if the handlebody deformation retracts onto the bouquet.  When the handlebody is equipped with a contact structure, we may additionally require the bouquet to be Legendrian.  Throughout, we will let $K_U$ and $K_V$ denote Legendrian skeletons of $U$ and $V$, respectively, and we let $\nu(K_U)$ and $\nu(K_V)$ denote their standard neighbourhoods.   

\begin{figure}[h]
\begin{center}
\includegraphics[scale=.8]{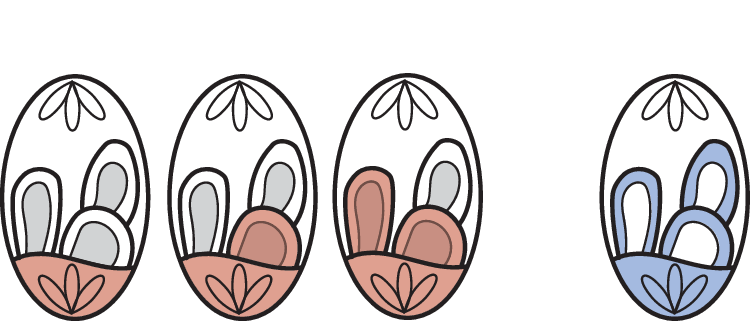}
\caption{Left: Adding successive bypasses produces a sequence of (red) handlbodies $U_i$ that are homeomorphic but not weakly contact isotopic to $\nu(K_U)$.  Right: Adding the  $1$-handles produces a bridge splitting $\widehat{H}=(\widehat{\Sigma}, \widehat{U}, \widehat{V})$; here $\widehat{U}$ is indicated in blue.}\label{fig:schematic}
\end{center}
\end{figure}

Begin with a smooth Heegaard splitting $\H=(\Sigma, U, V)$.  Since the complement $X=M\setminus \big(\nu(K_U)\cup\nu(K_V)\big)$ is topologically a product, Theorem \ref{thm:bypass} implies that up to weak contact isotopy, $X$ decomposes as a sequence of bypasses $X_1\ast \dots \ast X_f$ attached to $\partial \nu(K_U)$.\footnote{the triple $(\nu(K_U),\nu(K_V)),X)$ is called a "\emph{patty}" in \cite{BHH}. } 
Attaching only the first $i$ bypasses in the sequence determines a new handlebody $U_i$ with convex boundary; letting $i$ range over the indices of all bypasses produces a sequence of smoothly isotopic Heegaard surfaces $\Sigma_i:=  \partial U_i$.  Each $\Sigma_i$ has a bypass  arc $c_{i+1}$ that determines the bypass $X_{i+1}$.
Isotopy of $c_i$ through bypass arcs preserves the weak contact isotopy class of  $U_i$, so we may push each bypass arc off the faces of any $2$-handles already attached.  

\begin{definition} Let $X_1\ast\dots\ast X_f$ be a bypass decomposition of $M\setminus \big(\nu(K_U)\cup\nu(K_V)\big)$ as above.  Define $\widehat{U}$ to be the handlebody formed by adding only the bypass $1$-handles to $\nu(K_U)$:
\[\widehat{U}= \nu(K_U) \cup \big( \bigcup_{i=1}^f h^1_i \big).\] 
Then the \textit{bridge} of $\H$ is the Heegaard splitting
\[ \widehat{\H}( X_1\ast\dots\ast X_f)= \big(\partial \widehat{U}, \widehat{U}, M\setminus \text{int}(\widehat{U})\big).\]
\end{definition}

Attaching a contact $1$-handle preserves tightness, so it is immediate that $\widehat{U}$ is tight.  The complementary $\widehat{V}$ is also tight, as turning the manifold upside down builds $\widehat{V}$ from $1$-handles attached to $\nu(K_V)$.  In fact, the splitting $\widehat{\H}( X_1\ast \dots \ast X_f)=(\widehat{\Sigma}, \widehat{U}, \widehat{V})$ is convex, as co-core discs for each $1$-handle, together with meridional discs for the neighbourhoods of the skeleta, form a complete set of product discs.

\begin{remark}
Although the term ``bridge'' is new, this construction appears in both \cite{BHH} and \cite{LV} as a technique to construct a convex Heegaard splitting from an arbitrary smooth Heegaard splitting. 
\end{remark} 
When the original Heegaard surface $\Sigma$ is already convex in the splitting $\H=(\Sigma, U, V)$,  one may consider the relative bridge $\widehat{\H}^U$ built from a bypass decomposition of $M\setminus \big(\nu(K_U)\cup V\big)$, although this need not be convex.  Similarly, $\widehat{\H}^V$ is constructed from a bypass decomposition of $M\setminus \big(U \cup \nu(K_V)\big)$.  The relative bridges will be used in Proposition~\ref{prop:convexslice}.

\begin{proposition}\label{prop:extindep}
Up to positive stabilisation, the bridge \[\widehat{\H}(X_1\ast \dots \ast X_f)=(\widehat{\Sigma}, \widehat{U}, \widehat{V})\] is independent of the choice of bypass attaching arcs within their isotopy classes. 
\end{proposition}

\begin{proof} In order to construct a bridge splitting from a bypass decomposition, all the bypass attaching arcs must be isotoped off the faces of an $2$-handles already attached, but there is no canonical way to perform this isotopy. 
Although different isotopies preserve the weak contact isotopy type of each $U_i$, they may produce handlebodies which are inequivalent when the 2-handles are removed.  We will show that an isotopy of the attaching $0$-sphere for a $1$-handle on $\Sigma_i$ preserves the positive stabilisation class of  the bridge $\widehat{\H}$.  

\begin{figure}[h]
\begin{center}
\includegraphics[scale=1.2]{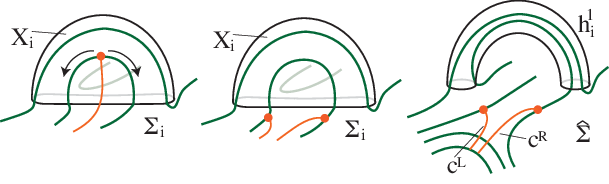}
\caption{ Left: Before constructing the bridge of $\H$, bypass arcs must be isotoped off the faces of the $2$-handles.  Centre: Isotoping $c$ across the face of a $2$-handle on $\Sigma_i$ preserves the weak contact isotopy class of $U_i$.  Right:  When $h_i^2$ is removed, the bypass arcs $c^R$ and $c^L$ are not isotopic on $\widehat{\Sigma}$.
} \label{fig:flowsphere}
\end{center}
\end{figure}

\begin{figure}[h]
\begin{center}
\includegraphics[scale=1.5]{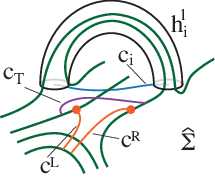}
\caption{ After attaching a bypass along $c_T$, the arcs $c^R$ and $c^L$ again become isotopic. } \label{fig:flowspherepf}
\end{center}
\end{figure}

Let $c^R$ and $c^L$ be two bypass attaching arcs for some $X_j$ ($j>i$) that are isotopic on $\Sigma_i$  via an isotopy that slides one endpoint across a face of some $2$-handle.  Because the proof is local, we may assume without loss of generality that the associated bypass is the most recent one; that is, $j=i+1$.  See Figure~\ref{fig:flowsphere}.  The Heegaard surfaces $\widehat{\Sigma}^R$ and $\widehat{\Sigma}^L$ associated to attaching bypasses along $c^R$ and $c^L$, respectively, are not convexly isotopic.  
 
 \begin{figure}[h]
\begin{center}
\includegraphics[scale=1.2]{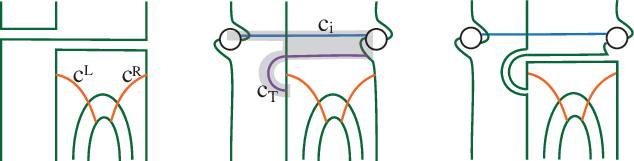}
\caption{ Left:  The bypass arcs $c^R$ and $c^L$ are isotopic on $\Sigma_i$.  Centre: If $h_i^2$ is not attached, then $c^R$ and $c^L$ are not isotopic. Right: After attaching a bypass along $c_T$, the arcs $c^R$ and $c^L$ are again isotopic.} 
\label{fig:inserttriv3big}
\end{center}
\end{figure}

Consider the bypass arc $c_T$ shown on Figure~\ref{fig:inserttriv3big}.  As the endpoints of $c^L$ and $c^R$ can be assumed to be arbitrarily close to $c_i$, we can also assume that $c_T$ is close to $c_i$.  Thus no other handles are attached on top of $c_T$ and it persists as a bypass arc on $\widehat{\Sigma}^L$. We claim that there is a  bypass half-disc attached along $c^T$ in $\widehat{V}^L$.

Observe that on $\Sigma_i$, the bypass arc $c_T$ is a rotation of $c_i$. As discussed in Section~\ref{sec:bypass}, there must exist a bypass disc $D_T$ attached along $c_T$.  Moreover, by Lemma~\ref{lem:rotex} $D_T$  may be assumed to lie in a neighbourhood of $X_i$ and the shaded area $R'$ in the centre picture of Figure~\ref{fig:inserttriv3big}.  It follows that the $1$-handle attached along $\partial c_L$ is disjoint from $D_T$, so $D_T$ persists in $\widehat{V}^L$.

Figure~\ref{fig:inserttriv3big} shows that  attaching the bypass along $c_T$ onto $\widehat{U}^L$ produces  a handlebody weakly contact isotopic to $\widehat{U}^R$. Then by Theorems~\ref{lem:obvshs} and \ref{thm:lv7.1}, the convex Heegaard decompositions $\widehat{\H}^L$ and $\widehat{\H}^R$ admit a common positve stabilisation. 
\end{proof}

Proposition~\ref{prop:extindep} implies that up to positive stabilisation,  the bridge splitting $\widehat{\H}=(\widehat{\Sigma},\widehat{U},\widehat{V})$   depends only on $K_U$, $K_V$, and the bypass decomposition $X=X_{1}\ast\cdots\ast X_{f}$.   As we show next, the bridge $\widehat{\H}$ is actually independent -- again, up to positive stabilisation -- of the bypass decomposition.

\begin{proposition}\label{prop:sliceindep}
    Given a smooth Heegaard splitting $\H=(\Sigma, U,V)$ let $K_U$ and $K_V$ be Legendrian skeletons for $U$ and $V$ and let $X=M\setminus \big(\nu(K_U)\cup\nu(K_V)\big)$.   Let $X_1\ast \dots \ast X_f$ and $X'_1\ast \dots \ast X'_{f'}$ be distinct bypass decompositions of $X$, and let $\widehat{\H}=\widehat{\H}(X_1\ast \dots \ast X_f) $ and $\widehat{\H}'=\widehat{\H}'(X'_1\ast \dots \ast X'_f) $ be the associated bridge Heegaard splittings formed by attaching  only the  $1$-handles of the respective bypasses  to $\nu(K_U)$. Then $\widehat{\H}$ and $\widehat{\H}'$ are related by a sequence of positive (de)stabilisations.
    \end{proposition}

\begin{proof}
According to Theorem~\ref{thm:bypass1p},  two bypass decompositions of $X$ are related by a finite sequence of isotopies, Trivial Insertion moves, and Far Commutation moves. Up to an isotopy that keeps the Heegaard surface convex,  Far Commutation preserves the handlebodies $\widehat{U}$ and  $\widehat{V}$. Thus we only need to show that inserting a trivial bypass preserves the positive stabilisation class of $\widehat{\H}$.  

Fix a bypass decomposition of $X$.
As above, let $\widehat{\H}=(\widehat{\Sigma},\widehat{U},\widehat{V})$ be the Heegaard splitting formed by attaching the associated $1$-handles to $\nu(K_U)$.  Now consider an alternative bypass decomposition that differs by the insertion of a trivial bypass $X_T=h^1_T\cup h^2_T$ attached along a trivial bypass arc $c_T$ somewhere in the sequence, and let $\widehat{\H}_T=(\widehat{\Sigma}_T,\widehat{U}_T,\widehat{V}_T)$ denote the associated Heegaard splitting.  Observe that the genus of $\widehat{\H}_T$ is larger than the genus of $\widehat{\H}$ because the additional $1$-handle $h_T^1$ is attached.

We claim that regardless of where the trivial bypass was inserted, $\widehat{\H}$ and $\widehat{\H}_T$ admit a common positive stabilisation. 
Recall that Proposition~\ref{prop:extindep} allows us to isotope each attaching sphere along respective components of  $\Gamma_{\Sigma_i}$ while preserving the positive stabilisation class of the bridge splitting. Using this flexibility, ensure that no handle attached after $h_T^1$ has an attaching sphere that intersects the disc supporting the trivial bypass.  Not only can any attaching sphere be isotoped off this region, but  we may also impose the stronger condition that subsequent attaching spheres are disjoint from the larger  disc shown in Figure~\ref{fig:isotdisc} that contains the bigon  certifying $c_T$ as a trivial bypass arc.  As a consequence, $c_T$ remains a trivial bypass arc on $\widehat{\Sigma}$;  we will use this below.   It follows that we may apply a sequence of Far Commutation moves to get a new bypass decomposition in which the bypass along $c_T$ is attached last. 

\begin{figure}[h]
\begin{center}
\includegraphics[scale=2]{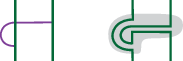}
\caption{ Left: The purple arc is a trivial bypass attachment arc.  Right: Attaching spheres for  handles attached after this bypass may always be isotoped off the shaded gray disc.} 
\label{fig:isotdisc}
\end{center}
\end{figure}

Attaching $h^1_T$  along $\partial c_T$ to $\widehat{U}$ produces the handlebody $\widehat{U}_T$. As noted in the comment before Theorem~\ref{thm:lv7.1}, adding a bypass $1$-handle to a convex splitting is always a positive stabilisation, which proves the stated claim.

\end{proof}

The following two results show that a convex Heegaard splitting $\H$ and its bridge $\widehat{\H}$ admit a common positive stabilisation.

\begin{proposition}\label{prop:convexslice}
    Given a tight Heegaard splitting $\H=(\Sigma, U,V)$, let $K_U$, $K_V$ be Legendrian skeletons for $U$ and $V$, respectively.   Let $X=X_{1}\ast\cdots \ast X_f$ be a bypass decomposition of the complement of $\nu(K_V)$ in $V$, where the index is the order in which the bypass slices are attached to $U$. Let $\widehat{\H}^U$ denote the relative bridge splitting with handlebodies $\widehat{U}=U\cup \big( \bigcup_{i=1}^f h^1_i \big)$ and $\widehat{V}=M\setminus \text{int}(\widehat{U})$.  Then the refinements of the Heegaard splittings $\widehat{\H}^U$ and $\H$ are related by a sequence of positive (de)stabilisations.
\end{proposition}

The relative bridge splitting $\widehat{\H}^U$ is tight by construction, so one may flip this Heegaard splitting upside down and apply Proposition~\ref{prop:convexslice} again to get the following consequence: 

\begin{corollary}\label{cor:cps}  Given a tight Heegaard splitting $\H=(\Sigma, U,V)$ let $K_U$, $K_V$ be Legendrian skeletons for $U$ and $V$ and let $X=X_{1}\ast\cdots\ast X_f$ be a bypass decomposition of the complement of $M\setminus\big(\nu(K_U) \cup \nu(K_V)\big)$.  Then the bridge splitting $\widehat{\H}$ and the refinement $\widetilde{\H}$  admit a common positive stabilisation.\end{corollary}

The proof of Proposition~\ref{prop:convexslice} is essentially the same argument as the proof of \cite[Theorem 6.1.]{LV}. Since it was not stated in full generality there, we indicate how it follows from the techniques of \cite{LV}.

\begin{proof}[Proof of Proposition~\ref{prop:convexslice}.] Define the sequence of intermediate Heegaard splittings $\widehat{\H}_j=(\widehat{\Sigma}_j,\widehat{U}_j,\widehat{V}_j)$ built by attaching only the first $j$ of the $1$-handles:
\[ \widehat{U}_j=U \cup \big( \bigcup_{i=1}^j h^1_i \big), \ \ \widehat{V}_j=M\setminus \text{int}(\widehat{U}_j), \ \  \widehat{\Sigma}_j=\partial \widehat{U}_j.\]

We will show that the refinements of successive Heegaard splittings in this sequence admit common positive stabilisations.  

Beginning with the case $j=1$, construct  the refinement  $\widetilde{\H}=(\widetilde{\Sigma}, \widetilde{U}, \widetilde{V})$ of $\H$. 

Recall that the refinement is constructed by drilling tunnels through meridional discs for the handlebodies.  We may choose these discs to be disjoint from the bypass $X_1=h^1_1\cup h^2_1$, so adding the $1$-handle $h_1^1$ to $U$ commutes with this drilling.  Thus attaching $h^1_1$ to the refined $\widetilde{U}$ produces the same Heegaard splitting as taking the refinement of $\widehat{\H}_1$.  As noted in Section~\ref{sec:ths}, attaching  $h_1^1$ in to $\widetilde{U}$ in the convex refinement is a positive stabilisation, which completes the base case.

Now proceed by induction on $j$.  For the inductive step, note that if $\widehat{\H}_j$ is tight, then $\widehat{\H}_{j+1}$ is also tight  because $\widehat{V}_{j+1}$ is a proper subset of the tight $\widehat{V}_j$, while $\widehat{U}_{j+1}$  is built from the tight $\widehat{U}_j$ by adding a contact $1$-handle.  The  argument above shows that refining and then adding a $1$-handle produces the same splitting as first adding the $1$-handle and then refining, which preserves the positive stabilisation class.  This establishes the inductive step.
\end{proof}

\subsection{Proof of the Giroux Correspondence}\label{sec:gc}
In order to prove the Giroux Correspondence, we consider a pair of convex Heegaard splittings of a fixed contact manifold.  We will modify these Heegaard splittings until they coincide, and each step in this evolution will preserve the positive stabilisation class of the Heegaard splittings.   

First, we repeat an argument from the proof of Theorem 6.1 in \cite{LV} that lets us restrict to  Heegaard splittings that are smoothly isotopic:

\begin{proposition}\label{lem:smoothstab}
   Any two convex Heegaard splittings  of $(M, \xi)$ have positive stabilisations  that are smoothly isotopic.
\end{proposition}

\begin{proof} In the topological setting, the Reidemeister-Singer Theorem asserts that any two Heegaard splittings of a fixed manifold become isotopic after sufficiently many stabilisations.  As there are no restrictions on the type of stabilisation that ensure this outcome, we are free to choose positive stabilisation of our convex Heegaard splittings. 
\end{proof}

Next, we can assume that the skeletons of the Heegaard splittings agree:

\begin{proposition}\label{lem:bouquet}
   Let $\H=(\Sigma,U,V)$ and $\H'=(\Sigma',U',V')$ be smoothly isotopic Heegaard splittings of $(M,\xi)$ with convex Heegaard surfaces. Then there exist Legendrian skeletons $K_U$ and $K_V$ and an isotopy $\Psi_t$ of $\H'$  keeping $\Psi_t(\Sigma')$ convex  such that $K_U$ is the skeleton of both $U$ and $\Psi_1(U')$ and $K_V$ is the skeleton of both $V$ and $\Psi_1(V)$.   
\end{proposition}

\begin{proof}

Choose Legendrian bouquets $K_*$ as the skeletons for each handlebody, where $*\in \{U,U',V,V'\}$.  By Theorem~\ref{thm:Legapprox1}, after possibly some stabilistaions we may assume that  the links $K_U\cup K_V$ and $K_{U'}\cup K_{V'}$ are Legendrian isotopic.


Extend this Legendrian isotopy to an ambient contact isotopy of the manifold, $\Psi_t$. By construction, for all $t$ the bouquet $\Psi_t(K_{{U'}})$ is a skeleton for the handlebody $\Psi_t({{U'}})$ with convex boundary $\Psi_t(\Sigma)$, and the analogous statement holds for $\Psi_t(K_{{V'}})$.  Since $\Psi_t$ is an extension of the Legendrian isotopy taking $K_{{U'}}\cup K_{{V'}}$ to $K_{{U}} \cup K_{{V}}$, we also have that $\Psi_1(K_{{U'}})=K_{{U}}$ is a skeleton for $U$ and $\Psi_1(K_{{V'}})=K_{{V}}$ is a skeleton for $V$, as desired. 
 \end{proof}

Finally, we are ready to prove the Giroux Correspondence for an arbitrary contact $3$-manifold.

\begin{proof}[Proof of Theorem \ref{thm:mainG}]
Let $\H$ and $\H'$ be two convex Heegaard splittings of $(M,\xi)$. By Proposition \ref{lem:smoothstab} we may assume that these Heegaard decompositions are smoothly isotopic. Then Proposition \ref{lem:bouquet} implies that after performing an isotopy keeping $\Sigma'$ convex, we may also assume there exist a Legendrian bouquet $K_U$  that is a  skeleton of both $U$ and $U'$ and a Legendrian bouquet $K_V$ that is a skeleton of both $V$ and $V'$. 

Choose bypass decompositions of the complement of the skeletons in $U$ and $V$:
\[{U\setminus \nu(K_U)} =
X_1\ast\cdots\ast X_l \text{ \ \ \ and \ \ \ } {V\setminus \nu(K_V)} =
X_{l+1}\ast\cdots\ast X_k.\] 

Concatenating these produces a bypass decomposition of $X=M\setminus \big(\nu(K_{U})\cup \nu(K_V)\big)$.  As usual, let $\widehat{U}$ denote the handlebody formed by attaching all of the $1$-handles to $\nu(K_U)$.  Setting $\widehat{\H}=(\partial \widehat{U}, \widehat{U}, M\setminus \text{int} (\widehat{U}))$, Corollary~\ref{cor:cps} implies that the refinement $\widetilde{\H}$ and the bridge $\widehat{\H}$ admit a common positive stabilisation.  To complete the argument, we recall Remark~\ref{rmk:extrastab}: since $\H$ is convex, any refinement of $\widetilde{\H}$ is a positive stabilisation of $\H$ itself.  Thus $\H$ and $\widehat{\H}$ admit a common positive stabilisation.

Now repeat this process using the alternative Heegaard splitting $\H'$:
 \[{U'\setminus \nu(K_U)}=
X'_1\ast\cdots\ast X'_{l'} \text{ \ \ \ and \ \ \ } V'\setminus \nu(K_V) =
X'_{l'+1}\ast\cdots\ast X'_{k'}.\]

Attaching all -- and only -- the $1$-handles to $\nu(K_U)$ gives a new bridge splitting $\widehat{\H}'$, and a similar argument shows that $\H'$ and $\widehat{\H}'$ admit a common positive stabilisation.  

However, $\widehat{\H}$ and $\widehat{\H}'$ are built from distinct bypass decompositions of $M\setminus \big(\nu(K_U) \cup \nu(K_V)\big)$, so Proposition~\ref{prop:sliceindep} implies that they are related by a sequence of positive (de)stabilisations, completing the proof. \end{proof}

\appendix  

\section{Bypass decompositions} \label{Appendix}


\hyperref[Appendix]{}
\subsection{Introduction}

This appendix presents a proof of Theorem~\ref{thm:bypass1p}, which is assumed in the main part of the paper. This result was first stated in the 3-dimensional case in \cite{BT}, while \cite{BHH} proves the analogous statement for all odd dimensions. Our approach  uses techniques from foliation theory which are specific to dimension three and in the spirit of Giroux \cite{Gi2,Gi3}.

Recall from Section~\ref{sec:bypassdecomp} that a contact structure $(\Sigma \times [0,1], \xi)$ with convex boundary admits a decomposition as $X_1*X_2*\dots*X_f$, where each $X_i$ is a topological product built by attaching a contact neighbourhood of a bypass half-disc to a thickening of the lower convex boundary.  Moves termed Far Commutation and Trivial Insertion in Section~\ref{sec:bypassdecomp} change this decomposition while preserving $(\Sigma \times [0,1], \xi)$ up to weak contactomorphism.  We restate Theorem~\ref{thm:bypass1p} for convenience:
\begin{theorem}\label{thm:main}[Theorem~\ref{thm:bypass1p}] Any two bypass decompositions of $(\Sigma\times I, \xi)$ are related by Bypass Arc Isotopy and a finite iteration of Trivial Insertion and Far Commutation moves. 
\end{theorem}

We give a high-level overview of the argument here, with the details deferred as indicated. See also Theorem \ref{thm:strongmain} for a stronger version of the result.

Section~\ref{sec:gnf} recalls the work of Giroux in \cite{Gi3} which defines a particularly nice type of a contact structure on $\Sigma\times [0,1]$.  He shows -- and we describe -- how to isotope an arbitrary contact structure to one in this \textit{Giroux normal form}. To such a contact structure, we associate a finite list of bypass arcs on selected convex surfaces $\Sigma \times\{t_i-\epsilon\}$.  Any such \textit{bypass sequence} (Definition~\ref{ref:bypseqfin}) determines a bypass decomposition in the sense of Theorem~\ref{thm:bypass}, and Definition~\ref{def:bseqequiv} declares two sequences to be equivalent if they are related by Far Commutation and Trivial Insertion.  We will argue that if $\xi$ and $\xi'$ are isotopic contact structures in Giroux normal form, then their bypass sequences are equivalent. 

The technical machinery propelling this work is the study of k-parameter families of surface foliations.   This is a natural arena to work in, as each slice $\Sigma\times \{t\}$ in the contact manifold inherits a characteristic foliation; $\Sigma\times [0,1]$ induces a 1-parameter family of such characteristic foliations; and an isotopy between contact structures on $\Sigma\times [0,1]$ gives rise to a 2-parameter family of characteristic foliations.  The topology of vector fields on a surface has been extensively studied, and we rely in particular on \cite{Sot} and \cite{KS}.  Importantly, certain spaces of vector fields are stratified, allowing us to characterise the kinds of foliations that may arise in generic k-parameter families.

The dynamics on planar surfaces are considerably simpler than on arbitrary surfaces; in order take advantage of this, our first intervention in Section~\ref{sec:sep} ``cuts''  an arbitrary $\Sigma$ into planar subsurfaces separated by annuli with prescribed characteristic foliations.  Although the surface itself can have any topology, imposing these local models forces the characteristic foliations to behave like foliations on planar surfaces. This construction, which generalises work of Giroux in \cite{Gi2}, makes the codimension 2 phenomena tractable. In the all-dimensions argument of \cite{BHH}, this outcome is achieved by inserting plugs to ensure the characteristic foliation is Morse-like at most times; this is the key difference between the two arguments.

Section~\ref{sec:strata} summarises essential results from the dynamical literature on the stratification of $C^k$ vector fields on planar surfaces.  

With the necessary tools assembled, Section~\ref{sec:gseq} explains how to associate an equivalence class of bypass sequences to any 1-parameter family of characteristic foliations that is both generic in the sense of \cite{KS} and separated in the sense of Section~\ref{sec:sep}.   This extends the assignment of a bypass sequence described earlier to a broader class of contact structures. In Section~\ref{sec:codim2} we consider a path of contact structures $\xi_s$ on $\Sigma \times I$.   Theorem~\ref{thm:equiv} asserts that the equivalence class of the bypass sequence is fixed along the path $\xi_s$.  The proof is a case chase: viewing the path $\xi_s$ as a 2-parameter family of foliations  $\mathcal{F}_{s,t}$,  the stratification described in Section~\ref{sec:strata} enumerates the characteristic foliations that can arise in the generic family $\mathcal{F}_{s,t}$ and we examine each in turn.  


Finally,  Section~\ref{sec:enhanced} introduces an \textit{enhanced bypass sequence} for a contact structure; this enriches the essentially combinatorial data of a bypass sequence with data about how the dividing set evolves with $t$ on $\Sigma\times[0,1]$.  The proofs may all be extended to this new object, proving the following result:

\begin{theorem}[Theorem~\ref{thm:strongmain}]

Suppose that $\xi$ and $\xi'$ are contact structures on $\Sigma \times I$  in Giroux normal form that are isotopic relative to $\Sigma \times \{0,1\}$. Then any enhanced bypass sequences admitted by $\xi$ and $\xi'$ are equivalent.  
\end{theorem}

Although this statement is structurally similar to Theorem~\ref{thm:equiv}, we note that enhanced bypass sequences are designed to distinguish contact manifolds that may produce identical bypass sequences.  This difference is explained in Section~\ref{sec:enhanced} and will be explored in forthcoming work.

This appendix assumes familiarity with the structure and properties of characteristic foliations, a topic developed by Giroux across many works: \cite{Gi,Gi2,Gi3}.  For background, we refer the non-Francophone reader to the exposition by Massot \cite{Ma} and to Mathews's translation \cite{Mat} of \cite{Gi}. We also make use of results and language from dynamics and suggest \cite{Sot, KS} as useful resources.

\subsection{Giroux normal form and bypass sequences}\label{sec:gnf}

Throughout this appendix, let $\xi$ denote a contact structure on $\Sigma\times [0,1]$ such that the boundary surfaces $\Sigma\times \{0,1\}$ are convex. 

\begin{definition}\label{def:Gnf} A contact structure $\xi$ on $\Sigma\times [0,1]$  is in \textit{Giroux normal form} if 
\begin{enumerate}
\item the movie of characteristic foliations $\mathcal{F}_t$ is a generic 1-parameter family;
\item  it is convex except at finitely many $t$-values $0<t_1<\dots<t_f<1$, and on each $\mathcal{F}_{t_i}$, the failure to be Morse-Smale is seen in a single retrograde connection. 
\end{enumerate}
\end{definition}

Observe that if $\xi$ is in Giroux normal form, then the  surfaces  $\Sigma\times \{0,1\}$ are convex. As explained next, we can associate a sequence of bypass arcs and dividing sets to any contact structure in Giroux normal form.  


\begin{definition}\label{def:bypasssequence} A \textit{strict bypass sequence} is a list $\mathcal{B}=(\Gamma_0, c_1, \Gamma_1, \dots, c_f, \Gamma_f)$ such that 
\begin{enumerate}
\item $\Gamma_i$ is a multicurve on $\Sigma$ that divides $\Sigma$ into two components;
\item  $c_i$ is a bypass arc on $(\Sigma,\Gamma_{i-1})$ with a distinguished disc neighbourhood $D_i$;
\item $\Gamma_i=\Gamma_{i-1}$ outside of $D_i$; and
\item within $D_i$, the multicurve $\Gamma_i$ is the result of altering $\Gamma_{i-1}$ by a bypass along $c_i$.  
\end{enumerate}
   \end{definition}

 \begin{figure}[h]
\begin{center}
\includegraphics[scale=1]{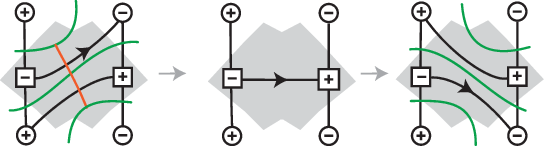} 
\caption{ Local model for the bypass arc associated to a retrograde connection.  Here, as elsewhere, $\Gamma$ is shown in green and the bypass arc is orange.  The disc $D$ is shaded. } 
\label{fig:retrogradeor}
\end{center}
\end{figure}

Figure~\ref{fig:retrogradeor} shows a standard local model for assigning a bypass arc to a neighbourhood of a retrograde orbit.  Using this, we may assign a strict bypass sequence to any contact structure in Giroux normal form. Set $\Gamma_0$ to be a dividing set on $\Sigma\times \{0\}$.  Since this changes only by isotopy in the interval $[0,t_1)$, pull the retrograde connection on $\Sigma\times \{t_1\}$ back by this isotopy 
and use the local model shown in Figure 1 to define $c_1$ (and consequently, $\Gamma_1$). Repeating this for each singular $t_i$-value produces a strict bypass sequence.

However, a strict bypass sequence is an overly rigid object, as it will suffice to consider dividing curves and bypass arcs up to isotopy. We clarify the allowable deformations next.

Given a strict bypass sequence, let $\Gamma_{f-1}$ deform by isotopy to $\Gamma'_{f-1}$. In addition, we may carry $c_f$ to  $c'_f$ via any isotopy through bypass arcs. Then define $\Gamma'_f$ as the image of $\Gamma'_{f-1}$ under a bypass along $c'_f$.  (It follows that $\Gamma'_f$ is isotopic as a multicurve to the original.)  
Inductively, we allow  isotopies of $\Gamma_{i-1}$ and isotopies of $c_i$ through bypass arcs, noting that these  deform $\Gamma_i$ and thus also 
the multicurve isotopy between $\Gamma_i$ and $\Gamma'_i$. 

\begin{definition}\label{ref:bypseqfin}A \textit{bypass sequence} is the equivalence class of a strict bypass sequence up to allowable deformations. \end{definition}

Although more cumbersome to define, a bypass sequence is more natural than a strict bypass sequence in the present setting.  For example, the construction described below Figure~\ref{fig:retrogradeor} produces a strict bypass sequence, but the isotopy between the dividing curves at different heights depends, a priori, on the choice of represenatives.  The space of dividing curves adapted to a given characteristic foliation is contractible, and the distinction between choices is absorbed into the isotopies permitted in a bypass sequence.  By definition, any representative of a bypass sequence may be deformed to a strict bypass sequence, so we lose no rigor by permitting the additional flexibility.



In the following, it will cause no confusion to  suppress the primes and simply use the notation introduced in  the definition of a  strict bypass sequence.  


Next, we present the operations of Far Commutation and Trivial Insertion in this new language.  
\begin{definition}
Suppose $D_i\cap D_{i+1}=\emptyset$. Then the strict bypass sequences
\[\dots ,\Gamma_{i-1},c_i, \Gamma_i, c_{i+1}, \Gamma_{i+1},\dots \text{ \ \ \ and \ \ \ } \dots, \Gamma_{i-1},c_{i+1}, \Gamma_i', c_{i}, \Gamma_{i+1},\dots\] are related by \textit{Far Commutation}.
\end{definition}

\begin{definition}
Suppose $c_T$ is a trivial bypass arc on $(\Sigma,\Gamma_i)$.  Set $\Gamma_T$ to be the multicurve produced by performing the bypass along $c_T$, and let  $\psi_T$ be an isotopy supported near the trivialising bigon that carries $\Gamma_T$ back to $\Gamma_i$.  Then the strict bypass sequences  \[\dots ,\Gamma_{i},c_i,\Gamma_{i+1},\dots \text{ \ \ \ and \ \ \ }  \dots, \Gamma_{i},c_T,\Gamma_{T},\psi_{T}c_i,\psi_{T}\Gamma_{i+1},\dots\] are related by \textit{Trivial Insertion}.
\end{definition}

Far Commutation and Trivial Insertion naturally descend to moves on bypass sequences: if two strict bypass sequences are related by one of these moves, then we declare their associated bypass sequences to be related by moves with the same names.  


\begin{definition}\label{def:bseqequiv} Two bypass sequences are \textit{equivalent} if they are related by a finite sequence of Far Commutations and Trivial Insertions.

\end{definition}

\begin{remark} \label{rmk:rot}
   As explained in Lemma~\ref{lem:rotex}, bypass rotation is a consequence of Trivial Insertion and Far Commutation.  Given $c_i$ in a bypass sequence $\mathcal{B}$, one may insert a rotation of $c_i$ immediately afterward while preserving the equivalence class of $\mathcal{B}$.
\end{remark}

\subsection{Separated foliations}\label{sec:sep}
The dynamics of foliations on spherical and planar surfaces are simpler than on surfaces with positive genus.  Specifically, genus zero obstructs sufficiently complicated bifurcations and limit sets.  In this section we explain how to modify the characteristic foliation on $\Sigma$ to decompose it into planar subsurfaces which may then be studied individually.   

After choosing a set of topological meridians that cut $\Sigma$ into planar subsurfaces, we modify the characteristic foliation in a neighbourhood of these curves to have parallel sets of \textit{separators} (Definition~\ref{def:sep}).  The separators ensure that as the foliations evolve, no trajectories can connect singularities on opposite sides of any meridian. For fixed $\xi$, this construction is Giroux's banalisation \cite[Lemma 2.10]{Gi2}. Proposition~\ref{prop:sep}  extends this idea  to a 1-parameter family $\xi_s$ that respects certain boundary conditions.   The result of this admittedly more complicated construction is an $(s,t)$-square of foliations with the property that at every interior point, the characteristic foliation on $\Sigma$ is  separated along some fixed set of meridians.  

\begin{definition}\label{def:sep} A \textit{separator} is a characteristic foliation on an embedded annulus $I\times S^1$ in which saddles and nodes are arranged in a $3\times 2n$ checkerboard configuration so that the sign is constant along each $3\times \{x\}$ arc spanning the annulus.   

\begin{figure}[h]
\begin{center}
\includegraphics[scale=.8]{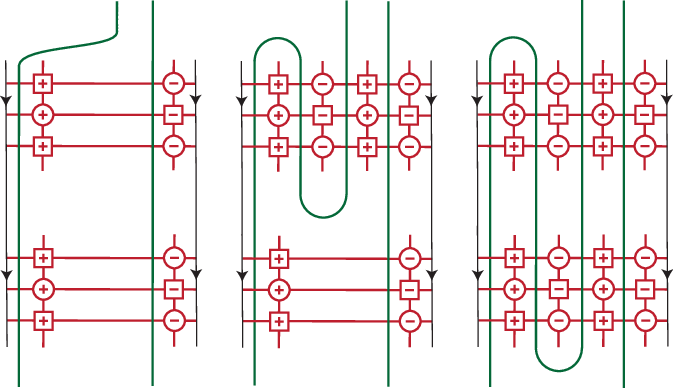}
\caption{Each annulus is a neighbourhood of a meridian $m$.  The left annulus has two parallel $3\times 2$ separators, while the right annulus has two parallel $3 \times 4$ separators.  The sequence shows three frames of a movie changing $|\Gamma\cap m|$.} 
\label{fig:sep}
\end{center}
\end{figure}
\end{definition}

See  Figure~\ref{fig:sep} for several examples of separators.

Let $M=\{m_1, \dots, m_m\}$ be a set of disjoint essential simple closed curves that cut the convex surface $\Sigma$ into planar subsurfaces.  We call $M$ a set of \textit{meridians}.

\begin{definition} The characteristic foliation on $\Sigma$ is \textit{separated} by $M$ if there exist one or more parallel separators in a product neighbourhood of each $m_i$.  A k-parameter family of characteristic foliations is \textit{separated} if there exists some $M$ such that each foliation is separated by $M$.
\end{definition} 

Given a convex surface $\Sigma$, it is always possible to arrange that the characteristic foliation is separated.  As we show next, this condition may also be imposed on $1$- and $2$-parameter families. Fix a set of meridians $M$.

\begin{proposition}\label{prop:sep} Suppose $\xi_s$ is a path of contact structures on $\Sigma\times [0,1]$ with the property that $\xi_0$ and $\xi_1$ are in Giroux normal form.   Then there exists an isotopic family $\xi'_s$ such that 
\begin{enumerate}
\item  $\xi_i'$  induces the same bypass sequence as $\xi_i$ for $i=0,1$;
\item away from a neighbourhood of $t=\{0,1\}$, the characteristic foliation on $\Sigma_{s,t}$ is separated by $M$; and 
\item on the neighbourhood of $t=\{0,1\}$, the surface $\Sigma_{s,t}$ is convex with respect to $\xi'_s$
.
\end{enumerate}
\end{proposition}

\begin{proof} It will be convenient to assume that $\xi_s$ is constant for $s$ near $\{0,1\}$, so we assume (or impose) this condition before proceeding in three steps following the schematic in Figure~\ref{fig:sep2}.  First, we normalise the contact structure near $t=0,1$, which is the region labeled $A$ in the figure.  We may assume that $\Sigma_{s,t}$ is convex for $t$ sufficiently close to $\{0,1\}$.  Within this region and closer to $t=\{0,1\}$, isotope $\mathcal{F}_{s,t}$ to be Morse-Smale and fix this throughout any later isotopies.  Isotope $\xi_s$, preserving the movies of bypass arcs associated to $\xi_0$ and $\xi_1$, so that for all  $(s,t)\in A$, the characteristic foliation on $\Sigma_{s,t}$ is separated by $M$.  More precisely, we require that near $t=1$, a neighbourhood of each $m_i$ has three parallel separators, while near $t=0$, a neighbourhood of each $m_i$ has two parallel separators.  This is doable by \cite[Remarque 2.7.]{Gi2},  
which tells us that the space of characteristic foliations adapted to a dividing set is contractible.  

\begin{figure}[h]
\begin{center}
\includegraphics[scale=.6]{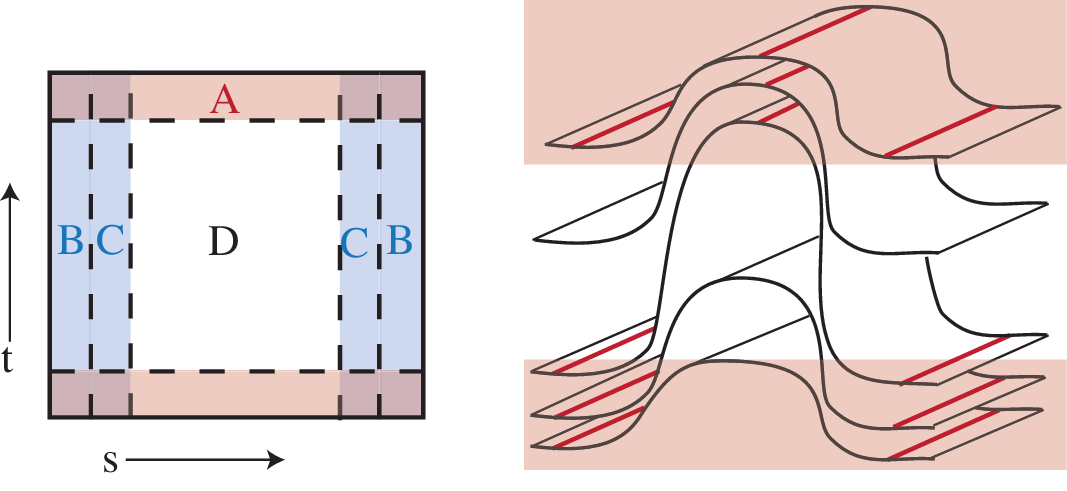}
\caption{Left: a schematic indicating modifications imposed on some initial 2-parameter family of foliations. Right: Vertical isotopy ensures that every of each $\Sigma\times \{t\}$ intersects at least one separator (red) in $A$.} 
\label{fig:sep2}
\end{center}
\end{figure}

Let $B\cup C$ be a subset of the region where $\xi_s$ is independent of $s$. Fixing $\xi'=\xi$ near $s=0,1$, we will adjust the characteristic foliation in $C$ so that it is separated by $M$ at every $t$ value. We work upwards from $A\cap C$, where this condition already holds.  If all bypass discs remain disjoint from a neighbourhood of $M$ as $t$ increases, then the characteristic foliation may be chosen to perpetuate the separators upward.  If some bypass disc intersects this neighbourhood, carefully constructed horizontal isotopies of the discs (described next) will render them disjoint; these preserve the height of the bypass discs, so the equivalence class of the  bypass sequence is preserved. Similarly, if an isotopy passes the dividing curve through $M$,  we apply the same horizontal isotopies.

Recall that a separator along $m_i$ has a $3\times 2k$ annular grid of singular points, where the value of $k$ counts the number of intersections between $\Gamma$ and $m_i$. Each $m_i$ has at least two parallel separators, so any change in $|\Gamma\cap M|$ may be realised by altering one separator at a time; this ensures that for all $t$, the characteristic foliations remains separated.  See, for example, Figure~\ref{fig:sep}. Now use $B$ to interpolate between $\xi_0$ and the foliations dictated on $C$.  

The final step of the modification is much more dramatic, and it may be helpful to consider the right-hand picture in  Figure~\ref{fig:sep2}.  Choose some value $t_0$ lying in the lower component of $A$ and isotope the neighbourhood of each $m_i$  on $\Sigma_{\frac{1}{2}, t_0}$ so that it is separated at $t=t_0$ by the original two separators but is separated at $t=1-t_0$ by the third separator.   Extend this to an isotopy defined on the rest of $\Sigma\times[0,1]$, keeping each surface $\Sigma_{\frac{1}{2},t}$ fixed away from a neighbourhood of $M$.  This isotopy ensures that with respect to $\xi_s$, every surface $\Sigma\times\{t\}$ is now separated by $M$. Finally, extend this isotopy to  all $(s,t)\in D$, using $s\in C$ to interpolate between the original ``flat'' surfaces and the newly mountainous ones.  The peaks for $\Sigma_{s,t}$ in $C$ may not reach up to $A$, but we have previously ensured that $\Sigma_{s,t}$ is separated by $M$ at every height $t$.  It follows that for $\Sigma_{s,t} \in C\cup D$, the characteristic foliation is separated by $M$.  
\end{proof}

The construction above may be modified so that $|\Gamma\cap m_i|$ is independent of $s$ and $t$.  Noting that $|\Gamma\cap m_i|=0\text{ mod }2$, it suffices to see that this value may be increased by two using only horizontal isotopy, as shown in Figure~\ref{fig:sep}.  Henceforth, we assume that this additional modification has been imposed.

After the modifications described above, reparametrise the square so that each point $(s,t)$ once again represents a specific surface.  Having standardised $\Gamma\cap m_i$ across all $\Sigma_{s,t}$ away from the boundary of the square, we may assume that the characteristic foliations on the boundary of some neighbourhood of $M$ are fixed for all  $(s,t)\in D$.  In the next section we consider modifying of a 2-parameter family of characteristic foliations, relative to the boundary of each complementary component of the neighbourhood of $M$.

\subsection{The stratification of $C^k$ vector fields}\label{sec:strata} 
This section describes a stratification on the space of characteristic foliations on a spherical or planar surface $\Sigma$.  We rely heavily on the work of \cite{KS} and \cite{Sot}, who study the more general case of arbitrary $C^k$ vector fields on planar surfaces.  Although they do not assume the non-isochore condition required of characteristic foliations, non-vanishing divergence is an open condition; thus it follows that for sufficiently small deformations, no isochore singularities may appear.  See also the discussion in \cite{Ma}. 

Following \cite{KS}, we assign an integer value ``cod'' to phenomena (singular points, connections) that appear in a foliation.   The \textit{expected codimension} of a characteristic foliation is then defined as the sum of the values associated to the constituent features with positive cod.  We appeal to the following result, recalled here only in the generality we require:

\begin{theorem}\cite[Theorems 1, 2]{KS}\label{thm:KS} In a generic $k$-parameter family of vector fields on the sphere with $k\leq 2$, let $F_n$ denote the parameters belonging to foliations with an expected codimension of $n$.  Then $F_n$ is an immersed submanifold of the parameter space with codimension $n$.  
\end{theorem}

The cases k=1 and k=2 will be examined in Sections~\ref{sec:gseq} and \ref{sec:codim2}, respectively. 
The analysis in \cite{KS} does not include closed orbits, but following \cite[Section 3.4]{Ma} , we briefly discuss them here for inclusion in the lists below. See \cite{Sot} for further details.  The degeneracy of a closed orbit is not discernible purely via  topology but rather, is defined in terms of the Poincar\'e first return map on a transversal to the orbit.   A closed orbit is \textit{generic} if the derivative of its Poincar\'e first  return map at 0 is not 1.  If the first derivative is $1$, then the orbit is degenerate.  Higher derivatives classify the degeneracy: an orbit is  $1$-degenerate if the second derivative is nonzero and it is 2-degenerate if the second derivative vanishes but the third does not.  A $k$-degenerate orbit is a cod=k phenomenon.

\begin{figure}[h]
\begin{center}
\includegraphics[scale=1]{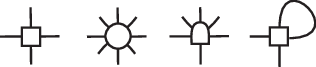}
\caption{Left: A saddle.  Centre left: A node.  Centre right: A saddle-node.  The separatrix pointing down is primary, while the horiztonal separatrices are secondary.  Right: A saddle loop, which is necessarily non-retrograde.} 
\label{fig:sing}
\end{center}
\end{figure}

The cod=0 phenomena are saddles, nodes, and nondegenerate orbits.  Since the codimension of a characteristic foliation is the sum of the cod values for each of its non-regular phenomena, foliations with only these limit sets and with no saddle-saddle connections are generic.  A generic characteristic foliation  ensures $\Sigma$ is convex.

Again following \cite[Table 1]{KS}, the features with cod=1 are as follows: 
\begin{enumerate}
\item\label{bif1} a  $1$-degenerate closed orbit; 
\item\label{bif2} a saddle-node; 
\item\label{bif3} a saddle-saddle connection; 
\item\label{bif4} a connection between a saddle and a primary or secondary separatrix of a saddle-node; 
\item\label{bif5}  a connection between primary or secondary separatrices of two saddle-nodes.
\end{enumerate}

It follows that a characteristic foliation in the codimension 1 stratum may have an unbounded number of saddles and nodes, along with exactly one of the phenomena from the list above.  In fact, the connections involving a saddle-node cannot appear in a codimension 1 foliation, as the existence of a saddle-node to anchor an end of this connection already contributes 1 to the count; thus, these cod=1 connections can only appear in a foliation in a stratum of higher codimension. 
\begin{figure}[h]
\begin{center}
\includegraphics[scale=.8]{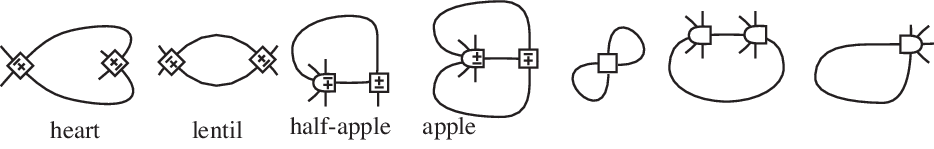}
\caption{Each polycycle above has codimension 2 because it exhibits two simultaneous cod=1 features.  Signs are indicated only in the named cases where  a retrograde connection is possible.} 
\label{fig:ensemb}
\end{center}
\end{figure}

Table 2 in \cite{KS} identifies the configurations with cod=2 as Bogdanov-Takens cusps, ultra-slow foci, degenerate saddles, and degenerate nodes.  The local models given there (Lemmas E and C in \cite{KS}) show that the first two of these are isochore and  thus will not arise in a characteristic foliation.  Since we will consider only strata of codimension less than or equal to 2, the remaining configurations may also be ignored: recall  Proposition 2.4 in \cite{Gi2} which states that a planar or spherical surface is convex if and only if its characteristic foliation has no degenerate orbits and no retrograde connections.  Since each of these obstructions to convexity contributes  cod$>$0, they cannot appear in a foliation that also has a cod=2 configuration.

Thus, relevant bifurcations in the codimension 2 stratum consist of either a 2-degenerate closed orbit or a pair of simultaneous cod=1 phenomena.  Enumerating the latter is a combinatorial task, which we summarise briefly.  

The less interesting cases involve cod=1 phenomena that are not combinatorially related, such as a retrograde connection and a saddle-node that doesn't connect to either of the endpoints.  This class includes foliations with a 1-degenerate orbit and another cod=1 feature, but the next section describes how to remove $1$-degenerate orbits.

The more interesting cases occur when there are combinatorial interactions between the cod=1 phenomena; the cases in which multiple singular points are cyclically ordered by their shared connections (i.e., \textit{polycycles}) were enumerated in \cite{KS} and are shown in Figure~\ref{fig:ensemb}.  The remaining cases are shown in Figure~\ref{fig:ensemb2}.

 \begin{figure}[h]
\begin{center}
\includegraphics[scale=.8]{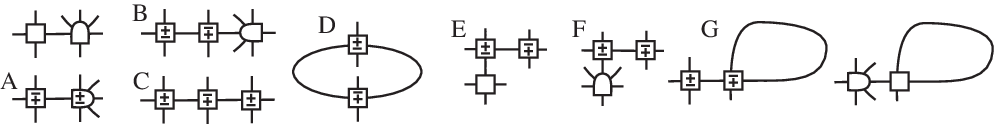}
\caption{Configurations that appear in the codimension 2 stratrum.  Signs are only indicated for configurations in which a retrograde connection is possible.}
\label{fig:ensemb2}
\end{center}
\end{figure}

We will use this classification to enumerate the bifurcations that appear in generic 1- and 2-dimensional families of characteristic foliations.

\subsection{Bypass sequences for separated characteristic foliations}\label{sec:gseq} 
Section~\ref{sec:gnf} explained how to assign a bypass sequence to a contact structure in Giroux normal form.  

\begin{proposition}\cite[Lemme~15]{Gi3}\label{prop:lemme15}  Let  $\xi$ be a contact structure on $\Sigma\times [0,1]$ such that the induced characteristic foliations on $\Sigma\times \{0,1\}$ are Morse-Smale. Then $\xi$ is isotopic relative to $\Sigma\times \{0,1\}$ to some $\xi'$  in Giroux normal form. 
\end{proposition}  

We may thus assign a bypass sequence to any contact structure with convex boundary, but it remains to show how this depends on the choice of isotopy.  Here, we recall the key feature of Giroux's proof of Proposition~\ref{prop:lemme15} in preparation for our discussion of invariance.  

According to the classification of codimension 1 phenomena in Section~\ref{sec:strata}, a generic 1-parameter family of characteristic foliations may have degenerate closed orbits, saddle-nodes, and saddle-saddle connections.  Saddle-nodes and non-retrograde connections are convex phenomena, and we have already seen how to assign a bypass arc to a retrograde connection.  This leaves degenerate closed orbits, which appear in movies where a pair of generic orbits cancel (death) or appear (birth), as shown in Figure~\ref{fig:degenor}. 

 \begin{figure}[h]
\begin{center}
\includegraphics[scale=.8]{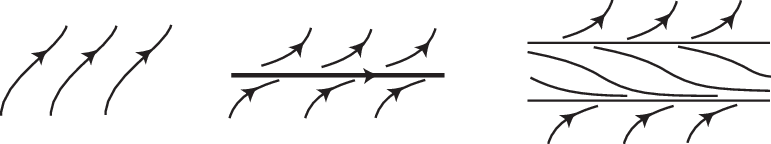}
\caption{ A movie showing the birth of two generic orbits via a degenerate orbit (bold). } 
\label{fig:degenor}
\end{center}
\end{figure} 

As in Giroux's Lemme 15 in \cite{Gi3}, suppose that the characteristic foliation has an isolated 1-degenerate orbit on $\Sigma\times \{t_0\}$.  One may deform the path of characteristic foliations on  $\Sigma\times (t_0-\epsilon, t_0+\epsilon)$ to obstruct the appearance of the degenerate orbit.  Specifically, introduce a pair of  saddle-nodes before $t_0$ and eliminate them after $t_0$, as shown in Figure~\ref{fig:frown1}.  Here, and elsewhere, we explicitly show only the case for a movie depicting the birth of two generic closed orbits; the death case may be seen by reversing the $t$-order and the orientation of the foliation.  

If the saddle-nodes persist throughout the modification, the configuration at $t_0$ is known as a \textit{lips bifurcation}, or more specifically, a \textit{malignant frown} \cite{KS} Figure 2b.  However, for $1$-parameter families of characteristic foliations it is possible to perturb the movie further, resolving each saddle-node to a saddle and a node as seen in the centre picture of Figure~\ref{fig:frown1}.   In this case the new path of foliations is generic and also convex away from  finitely many retrograde connections, as desired.  Furthermore, Giroux shows that the contact structure associated to this path of characteristic foliations is isotopic to the contact structure associated to the original path. 

 \begin{figure}[h]
\begin{center}
\includegraphics[scale=1]{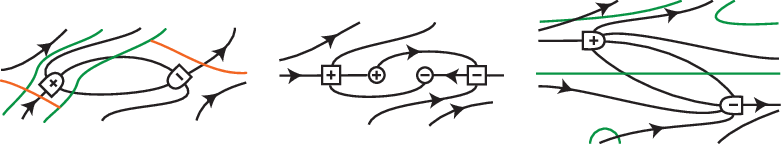}
\caption{ Inserting the malignant frown and resolving the saddle-nodes produces a generic path of foliations with no degenerate orbits (centre). 
} 
\label{fig:frown1}
\end{center}
\end{figure}

The frown insertion depends on a choice of $t$-interval supporting the deformation, and different choices will give rise to different bypass sequences.  For  the distinguished class of degenerate orbits we introduce next, Proposition~\ref{prop:tame} shows that the resulting bypass sequences are equivalent.   

In general, a degenerate closed orbit will be an accumulation point for a sequence of saddle-saddle connections.  This phenomenon was first identified in \cite{Sot} Remark 2.8.c, and the movie in Figure~\ref{fig:equivseq} indicates how such a sequence emerges as the orbits circulate along the annulus where the degenerate orbit will appear.  When the accumulating saddle-saddle connections are retrograde, each contributes a bypass arc to the bypass sequence, up  until the point when the malignant frown deformation is initiated.  Nevertheless, when the number of nearby saddle points is appropriately bounded, the equivalence classes of the bypass sequences associated to different choices agree.

\begin{figure}[h]
\begin{center}
\includegraphics[scale=1]{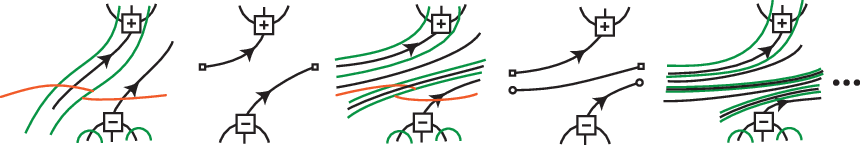}
\caption{ A sequence of retrograde connections limiting to a birth-type degenerate orbit.} 
\label{fig:equivseq}
\end{center}
\end{figure}

Given a degenerate orbit, let $m_\alpha$ denote the number of stable separatrices for positive saddles limiting away from the orbit and let $m_\omega$ denote the number of unstable separatrices for negative saddles limiting to the orbit.  

\begin{definition} A 1-degenerate closed orbit is \textit{tame} if $\min(m_\alpha, m_\omega)\leq 1$. 
\end{definition}

Now suppose that $\xi$ on $\Sigma\times[0,1]$ prints a characteristic foliation on $\Sigma_{t_0}$ that is Morse-Smale except for a 1-degenerate closed orbit.  Assume further that this is the unique degenerate orbit in the movie of characteristic foliations and that there are no retrograde orbits away from this annulus.   Choose $\epsilon>0$ such that the surfaces $\Sigma_{t_0\pm \epsilon}$ are convex and insert the malignant frown as in the construction above.   Thus modified, the movie has a finite number of retrograde connections supported in each of the three intervals $[0, t_0-\epsilon]$, $[t_0-\epsilon,t_0+\epsilon]$, and $[t_0+\epsilon,1]$.  Thus one can associate a bypass sequence to the modified movie. 

\begin{proposition}\label{prop:tame} For a tame orbit as above,  the bypass sequences associated to different choices of sufficiently small $\epsilon$ are equivalent. 
\end{proposition} 

When $\xi$ is in Giroux normal form, it was already possible to assign a bypass sequence.  Proposition~\ref{prop:tame} generalises the contact structures amenable to this process, at a cost of specifying only the equivalence class.  
\begin{proof}

Different choices of $\epsilon$ truncate the sequence of sparkling connections at different points. We will show that because the degenerate orbit is tame, the bypass arcs associated to any number of sparkling retrograde connections are rotations of the bypass arcs associated to the interval where the frown was inserted.  It follows that changing $\epsilon$ preserves the equivalence class of the  bypass sequence, as desired.  

When $\min(m_\alpha, m_\omega)=0$, no retrograde connections can occur, so $\epsilon$ is irrelevant. The more interesting case arises when $\min(m_\alpha, m_\omega)= 1$.  In this case $\max(m_\alpha, m_\omega)= 1$ counts the number of distinct pairs of saddles that can be connected by retrograde orbits, and the sparkling connections preserve a cyclic order on these pairs as orbits wind successively further around the annulus.  Figure~\ref{fig:equivseq} shows the bypass arcs associated to consecutive orbits in this sequence, while comparison to  Figure~\ref{fig:frown1} shows that these are all rotations of the bypass arc 
 associated to inserting the frown.  It follows from Remark ~\ref{rmk:rot} that appending any number of these arcs to the bypass sequence associated to the frown preserves its equivalence class, as desired.

The proofs for the other cases -- different signs, death -- are analogous. 

\end{proof}

\begin{remark} We believe that Proposition~\ref{prop:tame} could be proven at this point without the hypothesis that the degenerate orbit is tame. However, in the general case the sequence of sparkling connections is not unique, so the proof requires a more complicated combinatorial analysis. Since the stated version suffices for our purposes, we do not pursue this here. Of course, it is ultimately a consequence of Theorem \ref{thm:main}.
\end{remark}

\subsection{2-parameter families}\label{sec:codim2}

In this final section, we show that the equivalence class of the  bypass sequence assigned to a contact structure on $\Sigma\times[0,1]$ is independent of the choices made.

\begin{theorem}\label{thm:equiv} Suppose $\xi$ and $\xi'$ are contact structures on $\Sigma\times [0,1]$  that (1) are isotopic relative to the convex surfaces $\Sigma\times\{0,1\}$ and (2) induce a generic path of separated characteristic foliations on $\Sigma\times \{t\}$. 
Then any bypass sequences associated to $\xi$ and $\xi'$ are equivalent.
\end{theorem}

The previous section described how to remove a degenerate orbit by inserting a malignant frown, but recall that the resulting bypass sequence could depend on the $t$-interval supporting the modification.  Theorem~\ref{thm:equiv} establishes that this is not a source of concern, regardless of whether the original degenerate orbit was tame:  as the contact structures associated to different choices of $t$ are isotopic, their bypass sequences are equivalent.

We prove Theorem~\ref{thm:equiv} by assuming a 2-parameter family  $\mathcal{F}_{s,t}$ of characteristic foliations such that $\mathcal{F}_{s,0}$ is printed on $\Sigma$ by $\xi$ and $\mathcal{F}_{s,1}$ is printed on $\Sigma$ by $\xi'$.  After applying the construction in Section~\ref{sec:sep},  we assume that on the interior of the $(s,t)$ parameter space, all the foliations agree on the boundary of some neighbourhood of a set of meridians $M$.  We perturb the family $\mathcal{F}_{s,t}$ relative to these circles to get a generic family of foliations on each complementary planar subsurface. Next, we extend the construction in Section~\ref{sec:gseq}  to remove any non-tame $1$-degenerate closed orbits. These modifications can  be done with a perturbation preserving the isotopy class of the contact structures $\xi_t$. 

Away from the codimension 2 strata, the strict bypass sequences change only by isotopies, so the bypass sequence is preserved. Section~\ref{sec:strata}  enumerated the configurations appearing in foliations in the codimension 2 strata, and below, we verify that passing through any of these bifurcations preserves the equivalence class of the bypass sequence.  As many of the arguments proceed similarly, we provide a small number of them in detail and indicate how these provide a template for the remaining cases.
Working in a sufficiently small neighbourhood of the codimension 2 strata, we can assume that the dividing curve is unchanged away from the degenerate orbits and singularities.

In light of the discussion above, we will henceforth assume that $\mathcal{F}_{s,t}$ is a generic 2-parameter family of foliations on a planar surface,  fixed at the boundary. In each of the lemmata below, we propose that some non-generic phenomenon occurs at an interior point $(s_0, t_0)$ and examine the nearby foliations that arise.    Restricting to families of characteristic foliations reduces the number of cases; recall that by Giroux's crossing lemma \cite[Lemme 2.13.]{Gi2}, a movie with a retrograde connection corresponds to a bypass attachment only when it plays in the correct order, as shown in Figure~\ref{fig:retrogradeor}.  Each case shown generates parallel cases  by reflecting the foliation, reversing the orientation, and exchanging the signs of the singularities.

\begin{lemma}\label{lem:2codim1}  Suppose a 2-parameter family has a unique codimension 2 foliation at $(s_0,t_0)$ and furthermore, that $\mathcal{F}_{s_0,t_0}$ has  two combinatorially independent cod=1 phenomena.  Then the  bypass sequences associated to $\mathcal{F}_{s_0-\varepsilon,t}$ and $\mathcal{F}_{s_0+\varepsilon,t}$ are equivalent. 
\end{lemma}

\begin{proof}  Generically, two combinatorially independent cod=1 configurations appear in a single foliation when this is a frame in a family of movies exchanging the $t$-order in which they occur.  Because the two configurations are isolated from each other, any associated bypass arcs are necessarily disjoint, and hence the equivalence class of the bypass sequence is preserved if bypass arcs in the sequence occur in a different order. 
\end{proof}

\begin{lemma}\label{lem:IL}  Suppose a 2-parameter family has a unique codimension 2 foliation at $(s_0,t_0)$ and furthermore, that $\mathcal{F}_{s_0,t_0}$ has a configuration shown in Figure~\ref{fig:ensemb2}.  Then the bypass sequences associated to $\mathcal{F}_{s_0-\varepsilon,t}$ and $\mathcal{F}_{s_0+\varepsilon,t}$ are equivalent. \end{lemma}

\begin{proof}  
The non-labeled configurations in Figure~\ref{fig:ensemb2} have no retrograde connections or degenerate closed orbits;  the surface remains convex when these configurations appear, so the bypass sequence is necessarily preserved.

We examine each of the labeled foliations shown in Figure~\ref{fig:ensemb2} in turn. Configuration A occurs at the unique $s$ value when the primary separatrix of a saddle-node forms a retrograde connection to an existing saddle.  The distinct cases are enumerated by the number of saddles limiting to the node side of the saddle-node; Figure~\ref{fig:AB} shows this explicitly for two saddles.  The bypass sequence consists of  a single arc when the retrograde connection is between the original saddle and the one created by the saddle-node; this is the top path in Figure~\ref{fig:AB}.  Alternatively, each of the $k$ saddles that will eventually limit to the new node can separately form retrograde connections with the original saddle.  In this case, shown in the lower path, the bypass sequence consists of $k$ disjoint arcs.  The final one is isotopic to the arc associated to the other paths, and the preliminary arcs are rotations of it, so the  bypass sequences  are equivalent.  

The configurations labeled $B$ and $F$ each have only a single retrograde connection, and the associated bypass arc is independent of the synchronicity of the saddle-node.  

\begin{figure}[h]. 
\begin{center}
\includegraphics[scale=1]{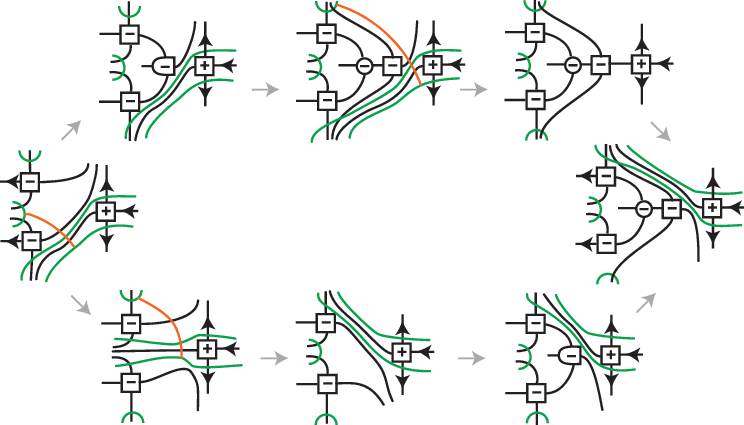}
\caption{ An indicative loop of characteristic foliations encircling the configuration labeled A in Figure~\ref{fig:ensemb2}.  The  bypass sequences associated to the two indicated paths are equivalent.} 
\label{fig:AB}
\end{center}
\end{figure}

Figure~\ref{fig:I} shows the equivalence of the bypass sequences for a family of movies passing through configuration C.  This configuration appears when the $t$-values of the two marked bypasses exchange as $s$ increases.  The bypass arcs are disjoint, so the equivalence class of the bypass sequences is independent of the order in which the retrograde connections occur, as desired.  Identifying the two positive saddles extends this argument to configuration D.

\begin{figure}[h]
\begin{center}
\includegraphics[scale=1]{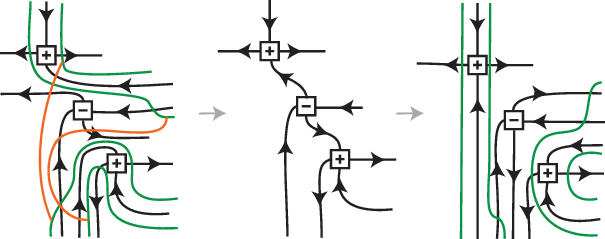}
\caption{ A movie of characterstic foliations passing through  configuration  C  Performing the indicated bypasses sequentially avoids the codimension 2 configuration.  Since the bypass arcs are disjoint, the  bypass sequences associated to performing these in either order are equivalent.} 
\label{fig:I}
\end{center}
\end{figure}

Configuration E breaks into two cases, depending on the sign of the saddle not involved in the retrograde connection, with one case shown in Figure~\ref{fig:E}.  Associated to the top row is a $\Gamma$ sequence with a single bypass arc, while the lower movie has two disjoint bypass arcs.  However, the Trivial Insertion turns the upper sequence into the lower sequence, so the  bypass sequences are equivalent.  

\begin{figure}[h]
\begin{center}
\includegraphics[scale=1]{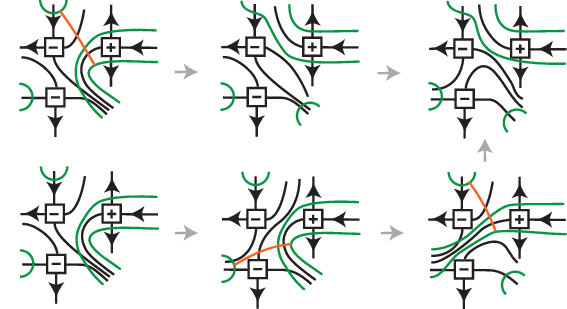}
\caption{ The two lefthand foliations are identical, and the arrows indicate distinct paths with equivalent  bypass sequences that lead to the upper right foliation.} 
\label{fig:E}
\end{center}
\end{figure}

A similar analysis for the $E$ configuration with alternating signs produces movies that each have a single (isotopic) bypass arc, so the associated  bypass sequences are equivalent.

\begin{figure}[h]
\begin{center}
\includegraphics[scale=1]{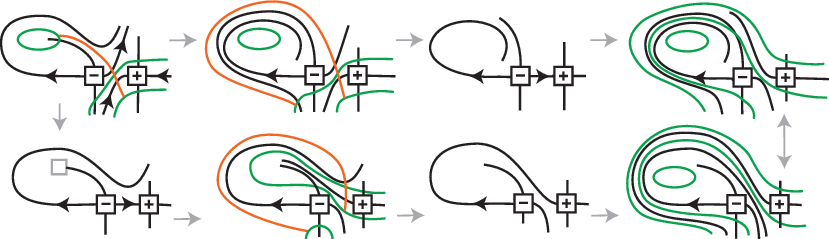}
\caption{ Two movies of characteristic foliations enclosing a snail.  The bypass arc indicated on the initial foliation in the upper left appears only in the bypass sequence associated to the lower path.  Since it represents a Trivial Insertion with respect to the bypass arc common to both paths,  the two  bypass sequences are equivalent. The grey shaded saddle in the lower right figure is discussed in the proof of Lemma~\ref{lem:acc1}.} 
\label{fig:snail}
\end{center}
\end{figure}

Finally, we examine the configuration labeled G in Figure~\ref{fig:ensemb2}, which we also refer to as the ``snail''.  Figure~\ref{fig:snail} shows that once again, the bypass sequences associated to two movies enclosing the snail differ only by Trivial Insertion and thus are equivalent. 
\end{proof}

The previous lemma addresses the case of a 1-parameter family of movies with a single codimension 2 configuration, but in general, these ensembles may be  accumulations point for other bifurcations.  When discussing bifurcations of a system of degenerate orbits and singularities below, we distinguish between  \textit{semi-local} and \textit{global} phenomena, where the former are completely contained within a small neighborhood of these cycles.

\begin{lemma}\label{lem:acc1} The bifurcations accumulating at foliations shown in Figure~\ref{fig:ensemb2} preserve the equivalence class of the bypass sequence. 
\end{lemma}

\begin{proof} 
The preceding proof examined loops of foliations enclosing the Figure~\ref{fig:ensemb2} configurations, and inspection reveals that in each case, there is a neighbourhood free from degenerate closed orbits.  Furthermore, the bifurcation diagrams already show all the interesting separatrices:  those whose  connections to the limit sets are parameter-dependent.  Any orbit or separatrix not completely contained within the shown neighborhood limits to/from a generic node or generic closed orbit and this connection persists for small changes in the parameter space. 

For each of the cases labeled A-F,  it follows that the bifurcation phenomena are particularly simple.  Direct inspection shows that an orbit leaving some subregion of the surface will not return to it, so no interesting accumulation phenomena arise. The snail (G) exhibits different behaviour, however, as an orbit limiting to the saddle loop will return to a neighbourhood of the saddle infinitely often.



The lower right picture in Figure~\ref{fig:snail} includes a shaded grey saddle in the interior of the loop of the snail.  If such a saddle exists, then a sequence of foliations with the codimension 2 configuration labeled E in Figure~\ref{fig:ensemb2} will accumulate at the snail.  However, Lemma~\ref{lem:IL} establishes that occurrences of this foliation preserve the equivalence class of the bypass sequence, as desired. 
\end{proof}

\begin{remark} We previously dismissed the occurrence of simultaneous combinatorially independent cod=1 phenomena as uninteresting.  Although this is true from the point of view of bypass sequences, it's worth briefly examining the case of a simultaneous saddle loop and retrograde connection from the dynamical point of view.  Codimension 2 foliations of type E may accumulate at this foliation.  This is not a source of concern, however, as  these accumulating foliations preserve the bypass sequence by Lemma~\ref{lem:IL}. 

\end{remark} 

Finally, we examine the polycycles shown in Figure~\ref{fig:ensemb}. 

\begin{lemma}\label{lem:heartetc}
Suppose that the generic family $\mathcal{F}_{s,t}$ contains a single configuration shown in Figure~\ref{fig:ensemb} or a single 2-degenerate orbit.  Then the bypass sequences associated to $\mathcal{F}_{s_0-\varepsilon,t}$ and $\mathcal{F}_{s_0+\varepsilon,t}$ are equivalent.
\end{lemma}

These bifurcations are distinguished by the fact that in each case, a neighbourhood in the parameter space contains a foliation with a closed 1-degenerate orbit.  It follows that a sequence of retrograde connections will accumulate at the named configuration, but in each case it is straightforward to verify that the degenerate orbits are tame, so there is no contradiction between these limit sets and the modifications previously made to remove non-tame orbits.  It then follows from Proposition~\ref{prop:tame}  that  a movie with one of these degenerate orbits admits a bypass sequence which is independent of when the malignant frown is inserted.

Each of the codimension 2 ensembles has been studied in depth elsewhere, and those authors provide bifurcation diagrams  describing the foliations that arise in a local neighbourhood in parameter space of the ensemble.  

\begin{proof}[Proof of Lemma~\ref{lem:heartetc}]  
We first consider  2-degenerate orbits.  Figure 3.13 in \cite{Ma} shows that a 2-degenerate orbit arises when two 1-degenerate orbits converge and cancel. Figure~\ref{fig:2degen} shows the movie before this cancellation has taken place, with the two 1-degenerate orbits in bold.  It is immediate that each of these 1-degenerate orbits is tame, so the movie admits a bypass sequence as shown. After the  cancellation, the movie is constant and hence, has no associated bypass arcs.  Thus, it suffices to show that the  bypass sequence shown in Figure~\ref{fig:2degen} is equivalent to the empty sequence.  Specifically, the bypass arc appearing in the middle frame is shown isotoped and shaded in the first frame.  Since the two bypass arcs admit disjoint  representatives, their order in the  bypass sequence may be exchanged via Far Commutation.  However, the shaded arc is produced via Trivial Insertion on the first frame, while the original first arc is a Trivial Insertion with respect to  the resulting dividing set.

\begin{figure}[h]
\begin{center}
\includegraphics[scale=.8]{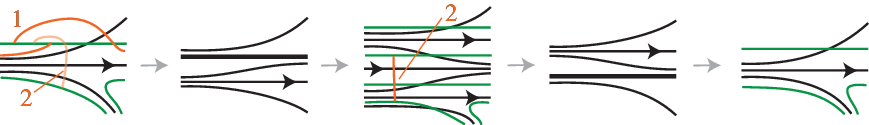}
\caption{ A movie of characteristic foliations with two tame 1-degenerate orbits. The associated bypass sequence is equivalent to the empty sequence.} 
\label{fig:2degen}
\end{center}
\end{figure}

The heart ensemble was studied by Dukov \cite{Duk}.  He distinguishes two cases, depending on whether the saddles have the same or different signs; we are concerned only with the latter, as the former generates no retrograde saddles.  Dukov shows that any neighbourhood of the heart contains a 1-degenerate orbit (or in his language, a ``parabolic cycle''), but we see that this closed orbit is tame: each of the saddles in the semi-local picture has a single separatrix limiting to this closed orbit, and orbits seen in the semi-local picture obstruct additional global saddles from limiting there. See Figure~\ref{fig:heart}.

The lower path is non-convex only at the indicated retrograde connection and hence the associated bypass sequence is the single arc shown.  The upper path has a sequence of sparkling retrograde connections accumulating at the degenerate orbit.   As this degenerate orbit is tame, it suffices by Proposition~\ref{prop:tame} to consider only the  bypass sequence associated to the lips.  Verifying that the  bypass sequences associated to the two paths are equivalent is then immediate. 

\begin{figure}[h]
\begin{center}
\includegraphics[scale=.7]{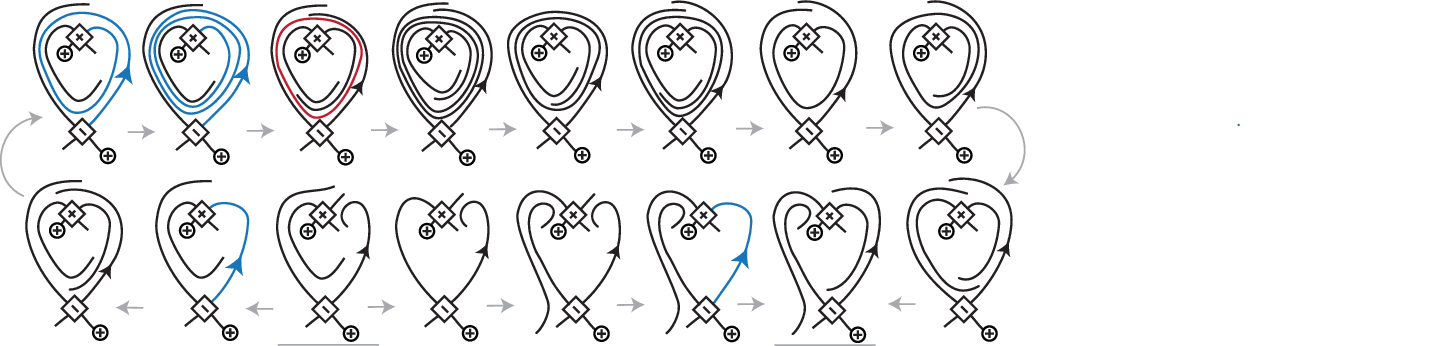}
\caption{ The two underlined  foliations are the starting and ending points for movies of characteristic foliations in a loop around the heart bifurcation.  The short path in the bottom row has a single retrograde connection (blue), while in the longer path, a family of sparkling retrograde connections (blue, indicative) accumulate at the degenerate orbit (red) } 
\label{fig:heart}
\end{center}
\end{figure}

In contrast to the heart, the lentil, apple, and half-apple admit relevant global phenomena.  We indicate how this affect the analysis of the half-apple and remark that other arguments are analogous.  

The black lines in Figure~\ref{fig:halfap} reproduce Figure 3 from  \cite{Gro} and show the semi-local foliations in a loop in parameter space that encloses the lentil.  The arrow identifies the degenerate orbit, and this is tame because the only saddle that limits to it from outside is the positive saddle already visible in the picture.  Expanding the frame of reference to include global connections introduces two additional phenomena, shown in grey.  First, a saddle outside the half-apple may have a separatrix that limits to the saddle-node.  For appropriate signs,  this may introduce retrograde connections, but these exist whether or not the half-apple occurs.  Second, there may be saddles inside the half-apple that limit to the saddle-node.  In this case, there will be a sequence of retrograde connections as the connecting separatrix spirals around inside the half-apple.  These will converge to the degenerate orbit, but as it is tame, Proposition~\ref{prop:tame} again applies.

Furthermore, the orientations of these additional separatrices preclude any interaction with each other, so only one separatrix of the positive saddle can connect to any of them; this obstructs any other codimension 2 phenomena from occurring.  Having established that no global connections materially alter the nearby foliations, we again assign a bypass sequence to each path around the half-apple and check that they are equivalent.

\begin{figure}[h]
\begin{center}
\includegraphics[scale=1]{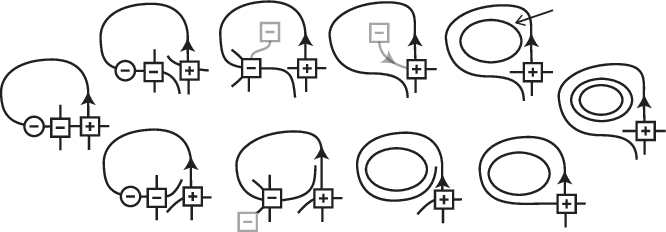}
\caption{ Foliations encircling the half-apple (\cite{Gro} Figure 3).} 
\label{fig:halfap}
\end{center}
\end{figure}

The semi-local bifurcations for the lentil  and apple are given in \cite{AAIS}, Figure 40a,  and \cite{Gro}, Figure 4, respectively, and they admit similar analysis.
\end{proof} 

\subsection{Enhanced bypass sequences}\label{sec:enhanced}

We have now proven Theorem~\ref{thm:main}, which is an essential ingredient in the proof of the Giroux Correspondence presented in the main part of the paper.  As we explain in this final section, a stronger version of the theorem may be established using essentially the same proof.

Observe that a non-contractible loop of isotopies of $\Sigma$ can be used to construct contactomorphic copies of $\Sigma\times [0,1]$ which are not contact isotopic relative to their boundary.  However, bypass sequences are blind to the distinction.  This poses no difficulty for the study of contact isotopic manifolds, but it limits the value of bypass sequences for more general applications.  Here, we upgrade the definition of a bypass sequence in order to detect such isotopies.  We similarly define an equivalence relation on these \textit{enhanced bypass sequences} to prove the following:

\begin{theorem}\label{thm:strongmain}
Suppose that $\xi$ and $\xi'$ are contact structures on $\Sigma \times I$  in Giroux normal form that are isotopic relative to $\Sigma \times \{0,1\}$. Then enhanced bypass sequences admitted by $\xi$ and $\xi'$ are equivalent.  
\end{theorem}

As we will show in a forthcoming paper, the equivalence classes of  enhanced bypass sequences are in bijection with contact structures in Giroux normal form, modulo isotopy fixing the boundary.  We introduce the notion of enhanced bypass sequences in the present work because the proof that this function is well defined is an easy generalisation of  the proof of Theorem~\ref{thm:main}.

A strict bypass sequence is a set of curves that change only by bypass attachments, while a bypass sequence allows each of these curves -- including the final ones -- to vary by isotopy.   In contrast, an \textit{enhanced bypass sequence} fixes the final dividing set.  This definition, and those that follow,  are the natural extensions of the objects introduced in Section \ref{sec:gnf} once this condition is imposed.

\begin{definition} An \textit{enhanced bypass sequence}  is a list \[\mathcal{EB} =\big(\Gamma_0(r), c_1, \Gamma_1(r),\dots, c_f, \Gamma_f(r)\big)\] where $\Gamma_i(r)$ is a $1$-parameter family of multicurves such that 
\begin{enumerate}
\item $\Gamma_i(r)$ divides $\Sigma$ into two components;
\item  $c_i$ is a bypass arc on $\big(\Sigma,\Gamma_{i-1}(1)\big)$ such that $c_i$ has a distinguished disc neighbourhood $D_i$;
\item $\Gamma_i(0)=\Gamma_{i-1}(1)$ outside of $D_i$; and
\item within $D_i$, the multicurve $\Gamma_i(0)$ is the result of altering $\Gamma_{i-1}(1)$ by a bypass along $c_i$. 
\end{enumerate}
In an enhanced bypass sequence,  $\Gamma_{start}:=\Gamma_0(0)$ and $\Gamma_{end}:=\Gamma_f(1)$.  
\end{definition}

As above, we  may assign an enhanced bypass sequence to a contact structure in Giroux normal form by watching how the dividing set changes away from the $t_i$-values where $\Sigma \times \{t_i\}$ is non-convex. This involves some choices, however, so we instead define a notion of compatibility between a contact structure and an enhanced bypass sequence.  

\begin{definition}\label{def:bseq} 

The contact structure $\xi$ on $\Sigma \times [0,1]$ in Giroux normal form \textit{admits the enhanced bypass sequence} $\mathcal{EB}$ if there exists $\epsilon>0$ such that the following hold:
\begin{enumerate}
\item $\Gamma_{start}$ divides the characteristic foliation on $\Sigma\times \{0\}$ and  $\Gamma_{end}$ divides the characteristic foliation on $\Sigma\times\{1\}$;
\item for $i=1,\dots,f$, there is an identification of the $i^{th}$ interval in the list \[[0,t_1-\epsilon]_t,[t_{i-1}+\epsilon, t_i-\epsilon]_t, [t_f+\epsilon, 1]_t\] with $[0,1]_r$ such that $\Gamma_{i-1}(r)$ divides $\Sigma\times \{r\}$; and 
\item\label{it:2} for $t\neq t_i$ in the interval $[t_i-\epsilon,t_i+\epsilon]$, the multicurve $\Gamma_i(0)\vert_{\Sigma\setminus D_i}$ can be extended as a dividing curve for the charactersitic foliation of $\Sigma\times \{t\}$.
\end{enumerate}
\end{definition}

\begin{definition}\label{def:ebseqequiv} Two enhanced bypass sequences $\mathcal{EB}$ and $\mathcal{EB}'$ are \textit{equivalent} if they are related by 
a finite sequence of the following: 
\begin{itemize} 
\item Isotopy

if $\mathcal{EB}$ and $\mathcal{EB}'$ are connected by a 1-parameter family of enhanced bypass sequences with $\Gamma_{start}$ and $\Gamma_{end}$ fixed throughout, then $\mathcal{EB}$ and $\mathcal{EB}'$ are equivalent; 
\item  Far Commutation

  if $D_i\cap D_{i+1}=\emptyset$ and if $\Gamma_i(r)$ and $\Gamma_i'(r)$ are $\ms{r}$-independent,  then 
\[\dots ,\Gamma_{i-1}(r),c_i, \Gamma_i, c_{i+1}, \Gamma_{i+1}(r),\dots \text{ \\ \ and \  \ } \dots, \Gamma_{i-1}(r),c_{i+1}, \Gamma_i', c_{i}, \Gamma_{i+1}(r),\dots\] are equivalent
; and
\item Trivial Insertion

 if $c_T$ is a trivial bypass arc on $(\Sigma,\Gamma_i)$  and $\Gamma_T(\ms{r})$ is an isotopy supported near $c_T$ and its trivialising bigon such that $\Gamma_T(1)=\Gamma_i(1)$,  then   \[\dots ,\Gamma_{i}(r),c_i,\Gamma_{i+1}(r),\dots \text{ \ \ \ and \ \ \ }  \dots ,\Gamma_{i}(r),c_T,\Gamma_{T}(r),c_i,\Gamma_{i+1}(r),\dots\] are equivalent.
\end{itemize}
\end{definition}

\begin{remark}  Note that an isotopy might carry the bypass arc $c_i$ through a disc $D_j$ for some $j<i$. 
\end{remark}

We close with a few remarks about the relationship between Theorem~\ref{thm:main} to Theorem~\ref{thm:strongmain}.  First, note that ignoring the data of $\Gamma_{end}$ defines  a forgetful map from equivalence classes of enhanced bypass sequences to equivalence classes of bypass sequences.  Applying this to Theorem~\ref{thm:strongmain} recovers  Theorem~\ref{thm:main}.

Although not required by the hypotheses of Theorem~\ref{thm:main},  at every step we allowed only contact isotopies relative to the boundary. In particular, these fix both the characteristic foliations $\mathcal{F}_0$ and $\mathcal{F}_1$. It follows that any representative of a bypass sequence used  in the main part of the proof is actually the image under the forgetful map of  a representative of an enhanced bypass sequence.  It is then straightforward to show that the equivalence class of bypass sequences associated to a 1-parameter family of contact structures in Section~\ref{sec:codim2} lifts to an equivalence class  of enhanced bypass sequences, as desired.

\bibliographystyle{alpha}
\bibliography{fob}
\end{document}